# A Kiefer–Wolfowitz theorem for convex densities

Fadoua Balabdaoui[1] and Jon A. Wellner[2,*]

*Université Paris-Dauphine and University of Washington*

**Abstract:** Kiefer and Wolfowitz [*Z. Wahrsch. Verw. Gebiete* **34** (1976) 73–85] showed that if $F$ is a strictly curved concave distribution function (corresponding to a strictly monotone density $f$), then the Maximum Likelihood Estimator $\widehat{F}_n$, which is, in fact, the least concave majorant of the empirical distribution function $\mathbb{F}_n$, differs from the empirical distribution function in the uniform norm by no more than a constant times $(n^{-1}\log n)^{2/3}$ almost surely. We review their result and give an updated version of their proof. We prove a comparable theorem for the class of distribution functions $F$ with convex decreasing densities $f$, but with the maximum likelihood estimator $\widehat{F}_n$ of $F$ replaced by the least squares estimator $\widetilde{F}_n$: if $X_1, \ldots, X_n$ are sampled from a distribution function $F$ with strictly convex density $f$, then the least squares estimator $\widetilde{F}_n$ of $F$ and the empirical distribution function $\mathbb{F}_n$ differ in the uniform norm by no more than a constant times $(n^{-1}\log n)^{3/5}$ almost surely. The proofs rely on bounds on the interpolation error for complete spline interpolation due to Hall [*J. Approximation Theory* **1** (1968) 209–218], Hall and Meyer [*J. Approximation Theory* **16** (1976) 105–122], building on earlier work by Birkhoff and de Boor [*J. Math. Mech.* **13** (1964) 827–835]. These results, which are crucial for the developments here, are all nicely summarized and exposited in de Boor [*A Practical Guide to Splines* (2001) Springer, New York].

## 1. Introduction: The Monotone Case

Suppose that $X_1, \ldots, X_n$ are i.i.d. with monotone decreasing density $f$ on $(0, \infty)$. Then the maximum likelihood estimator $\widehat{f}_n$ of $f$ is the well-known Grenander estimator: i.e. the left-derivative of the least concave majorant $\widehat{F}_n$ of the empirical distribution function $\mathbb{F}_n$.

In the context of estimating a decreasing density $f$ so that the corresponding distribution function $F$ is concave, Marshall [17] showed that $\widehat{F}_n$ satisfies $\|\widehat{F}_n - F\| \leq \|\mathbb{F}_n - F\|$ so that we automatically have $\sqrt{n}\|\widehat{F}_n - F\| \leq \sqrt{n}\|\mathbb{F}_n - F\| = O_p(1)$. Kiefer and Wolfowitz [14] sharpened this by proving the following theorem under strict monotonicity of $f$ (and consequent strict concavity of $F$). Let $\alpha_1(F) = \inf\{t : F(t) = 1\}$, and write $\|g\| = \sup_{0 \leq t \leq \alpha_1(F)} |g(t)|$.

**Theorem 1.1** (Kiefer–Wolfowitz [14]). *If $\alpha_1(F) < \infty$,*

$$\beta_1(F) \equiv \inf_{0 < t < \alpha_1(F)} (-f'(t)/f^2(t)) > 0,$$

*Supported in part by NSF Grant DMS-05-03822 and by NI-AID Grant 2R01 AI291968-04.
[1]Centre de Recherche, en Mathématiques de la Décision, Université Paris-Dauphine, Paris, France, e-mail: `fadoua@ceremade.dauphine.fr`
[2]University of Washington, Department of Statistics, Box 354322, Seattle, Washington 98195-4322, USA, e-mail: `jaw@stat.washington.edu`

$\gamma_1(F) \equiv \sup_{0<t<\alpha_1(F)} (-f'(t))/\inf_{0<t<\alpha_1(F)} f^2(t)) < \infty$, and $f'$ is continuous on $[0, \alpha_1(F)]$, then

$$\|\widehat{F}_n - \mathbb{F}_n\| = O((n^{-1} \log n)^{2/3}) \qquad \text{almost surely.} \tag{1}$$

Although Kiefer and Wolfowitz did not formulate their result in this way, the statement above follows from their proof. Also note that (1) implies that

$$\sqrt{n}\|\widehat{F}_n - \mathbb{F}_n\| = O(n^{-1/6}(\log n)^{2/3}) \to 0$$

almost surely, so that the MLE $\widehat{F}_n$ and the empirical distribution are asymptotically equivalent under the hypotheses of Theorem 1.

Kiefer and Wolfowitz [14] used Theorem 1.1 to show that the MLE $\widehat{F}_n$ of $F$ in the class of concave distributions is an asymptotically minimax estimator of $F$. (Also see Kiefer and Wolfowitz [15] for a generalization of the results of Kiefer and Wolfowitz [14] to allow somewhat weaker conditions.)

It follows from the rather general theorem of Millar [18] that the empirical distribution function $\mathbb{F}_n$ remains asymptotically minimax in a wide range of problems involving shape-constrained families of d.f.'s $\mathcal{F}$. In particular, for the classes $\mathcal{F}_k$ of distribution functions corresponding to $k$-monotone densities, it follows from Millar [18] that the empirical distribution function $\mathbb{F}_n$ is asymptotically minimax for estimation of $F$ even in the smaller classes $\mathcal{F}_k$. The interesting question which has not been addressed concerns asymptotic minimaxity of the MLEs within these classes. Our goal in this paper is to make some headway toward answering these questions by giving a partial (and imperfect) analogue of Theorem 1.1 in the case of $\mathcal{F}_2$, the class of distribution functions corresponding to the class of decreasing and convex densities. The MLE and least squares estimators of a density $f$ corresponding to $F \in \mathcal{F}_2$ have been studied by Groeneboom, Jongbloed and Wellner [11], and those results will provide an important starting point here.

In fact, we will not study the MLE, but its natural surrogate, the least squares estimator. This is because of the lack of a complete analogue of Marshall's lemma for the MLE in the convex case, while we do have such analogues for the least squares estimator; see Dümbgen, Rufibach and Wellner [7] and Balabdaoui and Rufibach [1].

One view of the Kiefer–Wolfowitz Theorem 1.1 is that it is driven by the (family of) corresponding local results, as follows:

**Theorem 1.2** (Local process convergence, monotone case). *Suppose that $t_0 \in (0, \infty)$ is fixed with $f(t_0) > 0$ and $f'(t_0) < 0$, and $f$ and $f'$ continuous in a neighborhood of $t_0$. Then*

$$n^{2/3}(\widehat{F}_n(t_0 + n^{-1/3}t) - \mathbb{F}_n(t_0 + n^{-1/3}t))$$
$$\Rightarrow \mathbb{C}_{b,c}(t) - \mathbb{Y}_1(t) \stackrel{d}{=} \left(\frac{2f^2(t_0)}{-f'(t_0)}\right)^{1/3} \{\mathbb{C}(at) - (W(at) - a^2t^2)\} \tag{2}$$

*in $(D[-K, K], \|\cdot\|)$ for every $K > 0$ where*

$$\mathbb{Y}_1(t) \equiv \sqrt{f(t_0)}W(t) + (1/2)f'(t_0)t^2 \equiv bW(t) - ct^2$$

*for $W$ a standard two-sided Brownian motion process starting from $0$, $\mathbb{C}_{b,c}$ is the Least Concave Majorant of $\mathbb{Y}_1$, $\mathbb{C} \equiv \mathbb{C}_{1,1}$ is the least concave majorant of $W(t) - t^2$, and $a \equiv \left([f'(t_0)]^2/(4f(t_0))\right)^{1/3}$.*

The (one-dimensional) special case of (2) with $t = 0$ is due to Wang [26], while the complete result is given by Kulikov and Lopuhaä [16].

Here the logarithmic term on the right side of (1) reflects the cost of transferring the family of (in distribution) local result to an (almost sure) global result. Here is a heuristic proof of (2); for the complete proof, see Kulikov and Lopuhaä [16]. For a similar result in the context of monotone regression, see Durot and Tocquet [8], and for a similar theorem in the context of the Wicksell problem studied by Groeneboom and Jongbloed [9], see Wang and Woodroofe [25]. For a related result in the context of estimation of an increasing failure rate, see Wang [24].

*Proof of Theorem 1.2.* We rewrite the left side of (2) as

$$
\begin{aligned}
& n^{2/3}\{\widehat{F}_n(t_0 + n^{-1/3}t) - \mathbb{F}_n(t_0 + n^{-1/3}t)\} \\
& \quad = n^{2/3}\{\widehat{F}_n(t_0 + n^{-1/3}t) - F(t_0) - n^{-1/3}f(t_0)t\} \\
& \quad\quad - n^{2/3}\{\mathbb{F}_n(t_0 + n^{-1/3}t)) - \mathbb{F}_n(t_0) - n^{-1/3}f(t_0)t\} \\
& \quad\quad + n^{2/3}\{\widehat{F}_n(t_0) - \mathbb{F}_n(t_0) - (\widehat{F}_n(\tau_0^-) - \mathbb{F}_n(\tau_0^-))\} \\
& \quad\quad - n^{2/3}\{\widehat{F}_n(t_0) - F(t_0)\}
\end{aligned}
\tag{3}
$$

where $\tau_0^-$ is the first point of touch of $\hat{F}_n$ and $\mathbb{F}_n$ to the left of $t_0$. From known local theory for $\hat{F}_n$ and $\mathbb{F}_n$ it follows easily that

$$
n^{2/3}\{\mathbb{F}_n(t_0 + n^{-1/3}t)) - \mathbb{F}_n(t_0) - n^{-1/3}f(t_0)t\}
$$
$$
\Rightarrow \sqrt{f(t_0)}W(t) + \frac{1}{2}f'(t_0)t^2 \equiv \mathbb{Y}_1(t), \tag{4}
$$

$$
n^{2/3}\{\widehat{F}_n(t_0 + n^{-1/3}t) - F(t_0) - n^{-1/3}f(t_0)t\} \Rightarrow \mathbb{C}_{b,c}(t) \tag{5}
$$

and

$$
n^{2/3}\{\widehat{F}_n(t_0) - F(t_0)\} \Rightarrow \mathbb{C}_{b,c}(0) \tag{6}
$$

where $\mathbb{C}_{b,c}$ is the least concave majorant of $\mathbb{Y}_1$. It remains to handle the third term. But since $\widehat{F}_n(t_0) - \widehat{F}_n(\tau_0^-) = \hat{f}_n(t_0)(t_0 - \tau_0^-)$ by linearity of $\widehat{F}_n$ on $(\tau_0^-, \tau_0^+)$,

$$
\begin{aligned}
& n^{2/3}\{\widehat{F}_n(t_0) - \mathbb{F}_n(t_0) - (\widehat{F}_n(\tau_0^-) - \mathbb{F}_n(\tau_0^-))\} \\
& \quad = -n^{2/3}(\mathbb{F}_n(t_0) - \mathbb{F}_n(\tau_0^-) - \hat{f}_n(t_0)(t_0 - \tau_0^-)) \\
& \quad = -n^{2/3}(\mathbb{F}_n(t_0) - \mathbb{F}_n(\tau_0^-) - f(t_0)(t_0 - \tau_0^-)) \\
& \quad\quad + n^{2/3}(\hat{f}_n(t_0) - f(t_0))(t_0 - \tau_0^-) \\
& \quad = n^{2/3}\{\mathbb{F}_n(t_0 + n^{-1/3}n^{1/3}(\tau_0^- - t_0)) - \mathbb{F}_n(t_0) \\
& \quad\quad - f(t_0)n^{-1/3}n^{1/3}(\tau_0^- - t_0)\} \\
& \quad\quad - n^{1/3}(\hat{f}_n(t_0) - f(t_0))n^{1/3}(\tau_0^- - t_0) \\
& \quad \to_d \mathbb{Y}_1(\tau_-) - \mathbb{C}_{b,c}^{(1)}(0)\tau_- = \mathbb{Y}_1(\tau_-) - \{\mathbb{C}_{b,c}(0) + \mathbb{C}_{b,c}^{(1)}(0)\tau_-\} + \mathbb{C}_{b,c}(0) \\
& \quad = \mathbb{Y}_1(\tau_-) - \mathbb{C}_{b,c}(\tau_-) + \mathbb{C}_{b,c}(0) = \mathbb{C}_{b,c}(0)
\end{aligned}
\tag{7}
$$

where $\tau_-$ is the first point of touch of $\mathbb{Y}_1$ and $\mathbb{C}_{b,c}$ to the left of 0, and hence $\mathbb{C}_{b,c}(\tau_-) = \mathbb{Y}_1(\tau_-)$. Combining (4), (5), (6) and (7) with (3) it follows that

$$
n^{2/3}\{\widehat{F}_n(t_0 + n^{-1/3}t) - \mathbb{F}_n(t_0 + n^{-1/3}t)\} \Rightarrow \mathbb{C}_{b,c}(t) - \mathbb{Y}_1(t)
$$

in $(D[-K, K], \|\cdot\|)$ for each fixed $K > 0$. □

## 2. The convex case

Now suppose that $X_1, \ldots, X_n$ are i.i.d. with monotone decreasing and convex density $f$ on $(0, \infty)$. Then the maximum likelihood estimator $\widehat{f}_n$ of $f$ is a piecewise linear, continuous and convex function with at most one change of slope between the order statistics of the data, and, as shown by Groeneboom, Jongbloed and Wellner [11], is characterized by

$$H_n(x, \widehat{f}_n) \begin{cases} \leq 1, \, x \geq 0 \\ = 1, \text{ if } \widehat{f}'_n(x-) < \widehat{f}'_n(x+). \end{cases}$$

where, with $\mathcal{K}$ being the class of convex and decreasing and nonnegative functions on $[0, \infty)$,

$$H_n(x, f) = \int_{[0,x]} \frac{2(x-y)/x^2}{f(y)} d\mathbb{F}_n(y), \qquad (x, f) \in \mathbb{R}^+ \times \mathcal{K}.$$

As shown by Groeneboom, Jongbloed and Wellner [11], the least squares estimator $\widetilde{f}_n$ of $f$ is also a piecewise linear, continuous, and convex function with at most one change of slope between the order statistics, but is characterized by

$$\widetilde{\mathbb{H}}_n(x) \begin{cases} \geq \mathbb{Y}_n(x), & x \geq 0, \\ = \mathbb{Y}_n(x), & \text{if } \widetilde{f}'_n(x-) < \widetilde{f}'_n(x+). \end{cases}$$

where $\widetilde{\mathbb{H}}_n(x) = \int_0^x \int_0^y \widetilde{f}_n(u) du dy \equiv \int_0^x \widetilde{F}(y) dy$ and $\mathbb{Y}_n(x) = \int_0^x \mathbb{F}_n(y) dy$. The corresponding estimators $\widehat{F}_n$ of $F$ and $Y$ are given by $\widehat{F}_n(x) = \int_0^x \widehat{f}_n(y) dy$ and $\widehat{\mathbb{H}}_n(x) = \int_0^x \widehat{F}_n(y) dy$ respectively. Since pointwise limit theory for both the MLE and the least squares estimators of $f$ are available from Groeneboom, Jongbloed and Wellner [11], we begin by formulating a (family of) local convergence theorems analogous to Theorem 1.2 in the monotone case. These will serve as a guide in formulating appropriate hypotheses in the context of our global theorem.

**Theorem 2.1** (Local process convergence, convex case). *If $f(t_0) > 0$, $f''(t_0) > 0$, and $f(t)$ and $f''(t)$ are continuous in a neighborhood of $t_0$, then for $(F_n, \mathbb{H}_n) = (\widehat{F}_n, \widehat{\mathbb{H}}_n)$ or for $(F_n, \mathbb{H}_n) = (\widetilde{F}_n, \widetilde{\mathbb{H}}_n)$,*

(8)
$$\begin{pmatrix} n^{3/5}(F_n(t_0 + n^{-1/5}t) - \mathbb{F}_n(t_0 + n^{-1/5}t)) \\ n^{4/5}(\mathbb{H}_n(t_0 + n^{-1/5}t) - \mathbb{Y}_n(t_0 + n^{-1/5}t)) \end{pmatrix}$$
$$\Rightarrow \begin{pmatrix} \mathbb{H}_2^{(1)}(t) - \mathbb{Y}_2^{(1)}(t) \\ \mathbb{H}_2(t) - \mathbb{Y}_2(t) \end{pmatrix}$$
$$\stackrel{d}{=} \begin{pmatrix} \left(24 \dfrac{f(t_0)^3}{f''(t_0)}\right)^{1/5} (\mathbb{H}_{2,s}^{(1)}(at) - \mathbb{Y}_{2,s}^{(1)}(at)) \\ \left(24^3 \dfrac{f(t_0)^4}{f''(t_0)^3}\right)^{1/5} (\mathbb{H}_{2,s}(at) - \mathbb{Y}_{2,s}(at)) \end{pmatrix}$$

*in $(D[-K, K], \|\cdot\|)$ for every $K > 0$ where*

$$\mathbb{Y}_2(t) \equiv \sqrt{f(t_0)} \int_0^t W(s) ds + \frac{1}{24} f''(t_0) t^4$$

and $\mathbb{H}_2$ is the "invelope" process corresponding to $\mathbb{Y}_2$: i.e. $\mathbb{H}_2$ satisfies: (a) $\mathbb{H}_2(t) \geq \mathbb{Y}_2(t)$ for all $t$; (b) $\int_{-\infty}^{\infty}(\mathbb{H}_2 - \mathbb{Y}_2)d\mathbb{H}_2^{(3)} = 0$; and (c) $\mathbb{H}_2^{(2)}$ is convex. Here

$$a = \left(\frac{f''(t_0)^2}{24^2 f(t_0)}\right)^{1/5},$$

and $\mathbb{H}_{2,s}$, $\mathbb{Y}_{2,s}$ denote the "standard" versions of $\mathbb{H}_2$ and $\mathbb{Y}_2$ with coefficients $1$: i.e. $\mathbb{Y}_{2,s}(t) = \int_0^t W(s)ds + t^4$.

Note that $\beta_2(F) \equiv \inf_{0 < t < \alpha_1(F)} (f''(t)/f^3(t))$ is invariant under scale changes of $F$, while $\delta_2(F) \equiv \sup_{0 < t < \alpha_1(F)} (f''(t)^2/f(t))^{1/5}$ is equivariant under scale changes of $F$; i.e. $\delta_2(F(c\cdot)) = c\delta_2(F)$.

*Proof.* Here is a sketch of the proof of the convergence in the first coordinate of (8). We write

$$n^{3/5}(F_n(t_0 + n^{-1/5}t) - \mathbb{F}_n(t_0 + n^{-1/5}t))$$
$$= n^{3/5}(F_n(t_0 + n^{-1/5}t) - F(t_0) - n^{-1/5}\frac{1}{6}f(t_0)t^3)$$
$$- n^{3/5}(\mathbb{F}_n(t_0 + n^{-1/5}t) - \mathbb{F}_n(t_0) - n^{-1/5}\frac{1}{6}f(t_0)t^3)$$
$$+ n^{3/5}(F_n(t_0) - \mathbb{F}_n(t_0) - (F_n(\tau_0^-) - \mathbb{F}_n(\tau_0^-)))$$
$$- n^{3/5}(F_n(t_0) - F(t_0)).$$

Here

$$n^{3/5}\left(F_n(t_0 + n^{-1/5}t) - F(t_0) - n^{-1/5}\frac{1}{6}f(t_0)t^3\right) \Rightarrow \mathbb{H}_2^{(1)}(t),$$
$$n^{3/5}\left(\mathbb{F}_n(t_0 + n^{-1/5}t) - \mathbb{F}_n(t_0) - n^{-1/5}\frac{1}{6}f(t_0)t^3\right) \Rightarrow \mathbb{Y}_2^{(1)}(t),$$
$$n^{3/5}(F_n(t_0) - F(t_0)) \Rightarrow \mathbb{H}_2^{(1)}(0),$$

while

$$n^{3/5}(F_n(t_0) - \mathbb{F}_n(t_0) - (F_n(\tau_0^-) - \mathbb{F}_n(\tau_0^-)))$$
$$= n^{3/5}\left(\mathbb{F}_n(t_0 + n^{-1/5}n^{1/5}(\tau_0^- - t_0)) - \mathbb{F}_n(t_0)\right.$$
$$\left. - n^{-1/5}\frac{1}{6}f(t_0)(n^{1/5}(\tau_0^- - t_0))^3\right)$$
$$- n^{3/5}\left(F_n(t_0 + n^{-1/5}n^{1/5}(\tau_0^- - t_0)) - F(t_0)\right.$$
$$\left. - n^{-1/5}\frac{1}{6}f(t_0)(n^{1/5}(\tau_0^- - t_0))^3\right)$$
$$+ n^{3/5}(F_n(t_0) - F(t_0))$$
$$\to_d \mathbb{Y}_2^{(1)}(\tau_-) - \mathbb{H}_2^{(1)}(\tau_-) + \mathbb{H}_2^{(1)}(0) = \mathbb{H}_2^{(1)}(0)$$

since $\mathbb{Y}_2^{(1)}(\tau_-) = \mathbb{H}_2^{(1)}(\tau_-)$. Combining the pieces yields the claim.
The proof for the second coordinate is similar. □

Now we can formulate our main result. Fix $\tau < \alpha_1(F)$. Our hypotheses are as follows:

R1. $F$ has continuous third derivative $F^{(3)}(t) = f''(t) > 0$ for $t \in [0, \tau]$ and $\beta_2(F, \tau) \equiv \inf_{0 < t < \tau}(f''(t)/f^3(t)) > 0$

R2. $\tilde{\gamma}_1(F, \tau) \equiv \sup_{0 < t < \tau}(-f'(t)/f^2(t)) < \infty$.

R3. $\gamma_2(F, \tau) \equiv \sup_{0 < t < \tau} f''(t) / \inf_{0 < t < \tau} f^3(t) < \infty$.

R4. $R \equiv \max\{1, \sup_{0 < t < \tau} f(t)\} / \inf_{0 < t < \tau} f(t) = \max\{1, f(0)\}/f(\tau) < \infty$.

In the rest of the paper we fix $\tau \in (0, \alpha_1(F))$ such that R1–R4 hold, and let $\|h\| \equiv \sup_{0 \le t \le \tau} |h(t)|$, the supremum norm of the real-valued function $h$ on $[0, \tau]$.

**Theorem 2.2.** *Suppose that R1–R4 hold. Then*

$$\|\widetilde{F}_n - \mathbb{F}_n\| \equiv \sup_{0 \le t \le \tau} |\widetilde{F}_n(t) - \mathbb{F}_n(t)| = O((n^{-1} \log n)^{3/5}), \tag{9}$$

$$\|\widetilde{\mathbb{H}}_n - \mathbb{Y}_n\| \equiv \sup_{0 \le t \le \tau} |\widetilde{\mathbb{H}}_n(t) - \mathbb{Y}_n(t)| = O((n^{-1} \log n)^{4/5}), \tag{10}$$

*almost surely.*

Note that (9) and (10) imply that

$$n^{1/2}\|\widetilde{F}_n - \mathbb{F}_n\| = O(n^{-1/10}(\log n)^{3/5}), \tag{11}$$

$$n^{1/2}\|\widetilde{\mathbb{H}}_n - \mathbb{Y}_n\| = O(n^{-3/10}(\log n)^{4/5}), \tag{12}$$

almost surely.

To prepare for the proof of Theorem 2.2, fix $0 < \tau < \alpha_1(F)$ for which the hypotheses of Theorem 2.2 hold. For an integer $k \ge 2$ define $a_j^{(k)} \equiv a_j \equiv F^{-1}((j/k)F(\tau))$ for $j = 1, \ldots, k$, and set $a_0^{(k)} \equiv a_0 \equiv \alpha_0(F) \equiv \sup\{x : F(x) = 0\}$. Note that $a_k^{(k)} = F^{-1}(F(\tau)) = \tau$ for all $k \ge 2$. We will often simply write $a_j$ for $a_j^{(k)}$, but the dependence of the knots $\{a_j\}$ on $k$ (and the choice of $k$ depending on $n$) will be crucial for our proofs. We also set $\Delta_j a = a_j - a_{j-1}$, and write $|a| = \max_{1 \le j \le k} \Delta_j a$.

Let $\mathbb{H}_{n,k}$ be the complete cubic spline interpolant of $\mathbb{Y}_n$ with knot points given by $\{a_j, j = 0, \ldots, k\}$. Thus $\mathbb{H}_{n,k}$ is piecewise cubic on $[a_{j-1}, a_j]$, $j = 1, \ldots, k$ with two continuous derivatives $\mathbb{H}_{n,k}^{(1)}$ and $\mathbb{H}_{n,k}^{(2)}$; see de Boor [5], pages 39–43 and 51–56. We will choose $k = k_n \sim (Cn/\log n)^{1/5} \to \infty$ in our arguments. $\mathbb{H}_{n,k_n}^{(2)}$ is not necessarily convex, but we will show that it becomes convex on $[0, \tau]$ with high probability as $n \to \infty$, and hence $\mathbb{H}_{n,k_n}$ will play a role analogous to the role played by the linear interpolation of $\mathbb{F}_n$ in the proofs of Kiefer and Wolfowitz [14]. (We will frequently suppress the dependence of $k = k_n$ on $n$, and write simply $k$ for $k_n$.)

Let $Y$ be defined by $Y(t) \equiv \int_0^t F(s) ds$; thus $Y^{(1)} = F$, $Y^{(j)} = f^{(j-2)}$, for $j \in \{2, 3, 4\}$. We will also need the complete cubic spline interpolant $H_{k_n}$ of $Y$; this will play the role of the linear interpolant $L = L^{(k)}$ of $F$ in Kiefer and Wolfowitz [14].

The cubic spline interpolant $\mathbb{H}_{n,k}$ of $\mathbb{Y}_n$ based on the knot points $\{a_j^{(k)}, j = 0, \ldots, k\}$ is completely determined on $[0, \tau]$ by the values of $\mathbb{Y}_n$ at the knots $a_j$, $j = 1, \ldots, k$ together with the values of $\mathbb{Y}_n^{(1)} = \mathbb{F}_n$ at 0 and $a_k = \tau$, namely $\mathbb{Y}_n(a_j)$, $j = 1, \ldots, J$, $\mathbb{Y}_n^{(1)}(0) = \mathbb{F}_n(0) = 0$, and $\mathbb{Y}_n^{(1)}(\tau)$; see, e.g., de Boor [5], page 43. As de Boor nicely explains in his Chapter IV, the complete cubic spline interpolant is one case of a family of cubic interpolation methods. Taking de Boor's function $g$ to be our present function $\mathbb{Y}_n$, several different piecewise cubic interpolants of $\mathbb{Y}_n$ can be described in terms of cubic polynomials $P_j$ on each of the intervals $[a_j, a_{j+1}]$ where the interpolating function $\mathbb{H}_n(\cdot; \underline{s})$ is given by $\mathbb{H}_n(x; \underline{s}) = P_j(x; \underline{s})$ for $x \in [a_j, a_{j+1}]$,

$j = 0, \ldots, k-1$, and where we require

$$P_j(a_j) = \mathbb{Y}_n(a_j), \qquad P_j(a_{j+1}) = \mathbb{Y}_n(a_{j+1})$$
$$P'_j(a_j) = s_j, \qquad P'_j(a_{j+1}) = s_{j+1},$$

for $j = 0, \ldots, k-1$. Here $\underline{s} = (s_0, \ldots, s_k)$ and the $s_j$'s are free parameters. Different choices of the $s_j$'s leads to different piecewise cubic functions agreeing with $\mathbb{Y}_n$ at the knots $a_j$; all of these different approximating functions $\mathbb{H}_n(\cdot; \underline{s})$ are continuous and have continuous first derivatives. Of interest to us here are the following particular ways of determining the $s_j$'s:

- $s_j = \mathbb{Y}_n^{(1)}(a_j) = \mathbb{F}_n(a_j)$, $j = 0, \ldots, k$. This gives the piecewise cubic *Hermite interpolant* of $\mathbb{Y}_n$, $\mathbb{H}_n(\cdot, \underline{s}) \equiv \mathbb{H}_{n,Herm}$.
- $s_j$, $j = 0, \ldots, k$ chosen so that $\mathbb{H}_n(\cdot, \underline{s}) \in C^2[0, \tau]$; i.e. so that $\mathbb{H}_n^{(2)}(\cdot, \underline{s})$ is continuous and $s_0 = \mathbb{Y}_n^{(1)}(0) = 0$ and $s_k = \mathbb{Y}_n^{(1)}(a_k) = \mathbb{Y}_n^{(1)}(\tau)$. This gives the *complete cubic spline interpolant* of $\mathbb{Y}_n$, $\mathbb{H}_n(\cdot, \underline{s}) \equiv \mathbb{H}_{n,CS} \equiv \mathbb{H}_{n,k}$.

The complete spline interpolant $\mathbb{H}_{n,CS}$ will play the role for us that the linear interpolant $L_n$ of $\mathbb{F}_n$ played in Kiefer and Wolfowitz [14]. As we will see, however, even though the Hermite interpolant $\mathbb{H}_{n,Herm}$ is not in $C^2[0, \tau]$ (i.e. $\mathbb{H}_{n,Herm}^{(2)}$ is not continuous), the slopes of its piecewise linear second derivative can be given explicitly in terms of $\mathbb{Y}_n$ and $\mathbb{Y}_n^{(1)} = \mathbb{F}_n$ at the knots, and our proof will proceed by relating the slopes of $\mathbb{H}_{n,Herm}^{(2)}$ to the (more complicated and less explicit) slopes of $\mathbb{H}_{n,CS}^{(2)} \equiv \mathbb{H}_{n,k_n}^{(2)}$ in order to prove point B in the following outline of our proof.

Here is an outline of the proof, paralleling the proof of the K–W theorem.

**Main steps, proof of (9) distribution function equivalence:**

**A.** By the generalization of Marshall's lemma for the convex density problem (see Dümbgen, Rufibach and Wellner [7]), for any function $h$ with convex derivative $h'$, $\|\widetilde{\mathbb{H}}_n^{(1)} - h\| \le 2\|\mathbb{F}_n - h\|$ where $\widetilde{\mathbb{H}}_n^{(2)} \equiv \widetilde{f}_n$. [This generalization is not yet available for the MLE $\widehat{\mathbb{H}}_n^{(1)}$ of $F$ in $\mathcal{F}_2$ corresponding to $\widehat{\mathbb{H}}_n^{(2)} = \widehat{f}_n$; see Dümbgen, Rufibach and Wellner [7] for a one-sided result.]

**B.** $P_F(A_n) \equiv P_F\{\mathbb{H}_{n,k_n}^{(2)}$ is convex on $[0, \tau]\} \nearrow 1$ as $n \to \infty$ if $k_n \equiv (C_0 \beta_2(F)^2 n / \log n)^{1/5}$ for some absolute constant $C_0$.

**C.** On the event $A_n$,

$$\begin{aligned}
\|\widetilde{\mathbb{H}}_n^{(1)} - \mathbb{F}_n\| &= \|\widetilde{\mathbb{H}}_n^{(1)} - \mathbb{H}_{n,k_n}^{(1)} + \mathbb{H}_{n,k_n}^{(1)} - \mathbb{F}_n\| \\
&\le 2\|\mathbb{F}_n - \mathbb{H}_{n,k_n}^{(1)}\| + \|\mathbb{H}_{n,k_n}^{(1)} - \mathbb{F}_n\| \\
&\qquad \text{by the generalization of Marshall's lemma (A)} \\
&= 3\|\mathbb{F}_n - \mathbb{H}_{n,k_n}^{(1)}\| \\
&= 3\|\mathbb{F}_n - \mathbb{H}_{n,k_n}^{(1)} - (F - H_{k_n}^{(1)}) + F - H_{k_n}^{(1)}\| \\
&\le 3\|\mathbb{F}_n - \mathbb{H}_{n,k_n}^{(1)} - (F - H_{k_n}^{(1)})\| + 3\|F - H_{k_n}^{(1)}\| \\
&\equiv 3D_n + 3E_n.
\end{aligned}$$

**D.** We show that $D_n = O((n^{-1} \log n)^{3/5})$ almost surely via a generalization of the K–W Lemma 2. We also show that $E_n = O((n^{-1} \log n)^{3/5})$ by an analytic (deterministic) argument.

Of course proving step **B** in this outline involves showing that the slopes of the $\mathbb{H}_{n,k_n}^{(2)}$ become ordered with high probability for large $n$, and this explains our interest in the slopes of both $\mathbb{H}_{n,CS}^{(2)} = \mathbb{H}_{n,k_n}^{(2)}$ and $\mathbb{H}_{n,Herm}^{(2)}$.

The assertion (10) of Theorem 2.2 can be proved in a similar way if we replace $\widetilde{\mathbb{H}}_n^{(1)}$, $\mathbb{H}_{n,k_n}^{(1)}$, $H_{k_n}^{(1)}$, $\mathbb{F}_n$, $F$ by $\widetilde{\mathbb{H}}_n$, $\mathbb{H}_{n,k_n}$, $H_{k_n}$, $\mathbb{Y}_n$, $Y$ respectively, and if we replace **A** by the following recent result of Balabdaoui and Rufibach [1]:

**A′**. For any function $G$ with convex second derivative $g''$, $\|\widetilde{\mathbb{H}}_n - G\| \leq \|\mathbb{Y}_n - G\|$.

*Proof of (9) assuming B.* First the deterministic term $E_n$. As in de Boor [5], page 43, let $I_4$ denote the complete cubic spline interpolation operator, and (as in de Boor [5], page 31, let $I_2$ be the piecewise linear (or "broken line") interpolation operator. Then by de Boor [5], (20) on page 56, with $p_n \equiv 1/k_n$,

$$E_n = \|F - H_{k_n}^{(1)}\| = \|Y^{(1)} - (I_4 Y)^{(1)}\| \leq \frac{1}{24}|a|^3 \|Y^{(4)}\|$$
$$\leq \frac{1}{24}\gamma_2(F,\tau)p_n^3 = O((n^{-1}\log n)^{3/5}).$$

To handle $D_n$, let $\$_3$ be defined to be the space of all quadratic splines on $[0,\tau]$, and similarly let $\$_2$ be the space of all linear splines on $[0,\tau]$. Then, by de Boor [5], page 56, equation (17), together with (18) on page 36, it follows that with

$$D_n = \|\mathbb{F}_n - \mathbb{H}_{n,k_n}^{(1)} - (F - H_{k_n}^{(1)})\| = \|(\mathbb{Y}_n - Y)^{(1)} - (I_4(\mathbb{Y}_n - Y))^{(1)}\|$$
$$\leq \frac{19}{4}\text{dist}((\mathbb{Y}_n - Y)^{(1)}; \$_3) \leq \frac{19}{4}\text{dist}((\mathbb{Y}_n - Y)^{(1)}; \$_2)$$
$$\leq \frac{19}{4}\|(\mathbb{Y}_n - Y)^{(1)} - I_2[(\mathbb{Y}_n - Y)^{(1)}]\|$$
$$= \frac{19}{4}\|(\mathbb{F}_n - F) - I_2(\mathbb{F}_n - F)\|$$
$$\leq \frac{19}{4}\omega(\mathbb{F}_n - F; |a|)$$
$$\stackrel{d}{=} \frac{19}{4}n^{-1/2}\omega(\mathbb{U}_n; p_n)$$
$$= O(n^{-1/2}\sqrt{p_n \log(1/p_n)}) \qquad \text{almost surely}$$
$$= O((n^{-1}\log n)^{3/5});$$

here

$$\omega(g;h) \equiv \sup\{|g(t) - g(s)| : |t - s| \leq h\},$$
$$\text{dist}(g;S) \equiv \min\{\|g - f\| : f \in S\}, S \subset C[0,\tau] \quad \text{and}$$
$$\mathbb{U}_n(t) \equiv \sqrt{n}(\mathbb{G}_n(t) - t)$$

where $\mathbb{G}_n(t) = n^{-1}\sum_{i=1}^n 1_{[0,t]}(\xi_i)$ is the empirical distribution function of $\xi_1, \ldots, \xi_n$ i.i.d. Uniform(0, 1) random variables. (See de Boor [5], pages xviii, 24, 32, and 34 for definition and use of $\text{dist}(g;S)$ and the modulus of continuity $\omega$ in conjunction.) □

*Proof of (10) assuming B.* By Hall [12] (also see Hall and Meyer [13] for optimality of the constant and de Boor [5], page 55),

$$E_n \equiv \|Y - H_{k_n}\| \leq \frac{5}{384}|a|^4\|Y^{(4)}\| \leq \frac{5}{384}R\gamma_2(F)\frac{1}{k_n^4}.$$

To handle the first term $D_n$, we note that

$$\mathbb{Y}_n - Y - (\mathbb{H}_{n,k_n} - H_{k_n}) = (\mathbb{Y}_n - Y) - I_4(\mathbb{Y}_n - Y)$$

where $I_4$ is the complete spline interpolant, and, on the other hand, for any differentiable function $g$ it follows from de Boor [5], page 45, equation (14), together with (18) on page 36, that

$$\begin{aligned}\|g - I_4 g\| &\leq \frac{19}{8}|a|\mathrm{dist}(g', \$_3) \leq \frac{19}{8}|a|\mathrm{dist}(g', \$_2) \\ &\leq \frac{19}{8}|a|\|g' - I_2 g\| \leq \frac{19}{8}|a|\omega(g', |a|).\end{aligned}$$

Applying this to $g = \mathbb{Y}_n - Y$, it follows that

$$\begin{aligned}\|\mathbb{Y}_n - Y - (\mathbb{H}_{n,k_n} - H_{k_n})\| &= \|(\mathbb{Y}_n - Y) - I_4(\mathbb{Y}_n - Y)\| \\ &\leq \frac{19}{8}|a|\,\omega(\mathbb{F}_n - F, |a|) \\ &\stackrel{d}{=} n^{-1/2}\omega(\mathbb{U}_n, p_n)\end{aligned}$$

Therefore $\omega(\mathbb{F}_n - F; |a|) = O(n^{-1/2}\sqrt{p_n \log(1/p_n)})$ almost surely (just as in the proof of Lemma 2 for the Kiefer-Wolfowitz theorem, see Section 5), we see that the order of $D_n$ is

$$n^{-1/2} p_n^{3/2} (\log(1/p_n))^{1/2} = O((n^{-1} \log n)^{4/5}) \qquad \text{almost surely}$$

as claimed. Thus the claim (10) is proved if we can verify that **B** holds. □

We end this section with a short list of further problems:

- It would be of interest to prove a comparable theorem for the MLE $\widehat{F}_n$ itself rather than $\widetilde{F}_n$. This involves several additional challenges, among which is a complete analogue of Marshall's lemma.
- Are either $\widetilde{F}_n$ or $\widehat{F}_n$ asymptotically minimax for estimating $F \in \mathcal{F}_2$?
- We conjecture that similar results hold for $k$−monotone densities and corresponding distribution functions ($k = 1$ corresponds to the Kiefer and Wolfowitz monotone density case, while $k = 2$ corresponds to the convex density case treated here). More concretely, we conjecture that under comparable hypotheses

$$\|F_n - \mathbb{F}_n\| = O((n^{-1} \log n)^{(k+1)/(2k+1)}) \qquad \text{almost surely}$$

  for $F_n = \widetilde{F}_n$ or $F_n = \widehat{F}_n$, the least squares estimator or MLE of $F \in \mathcal{F}_k$. Some progress on the local theory of the corresponding density estimators is given in Balabdaoui and Wellner [2] and Balabdaoui and Wellner [3]. On the interpolation theory side, the results of Dubeau and Savoie [6] may be useful.
- What is the exact order (in probability or expectation) of $\|\widehat{F}_n - \mathbb{F}_n\|$ in the case $k = 2$? Is it $(n^{-1} \log n)^{3/5}$ as perhaps suggested by the results of Durot and Tocquet [8] in the case $k = 1$?

## 3. Asymptotic convexity of $\mathbb{H}^{(2)}_{n,k_n}$

In this section we write $\mathcal{C}$ for the complete spline interpolation operator that maps functions $g \in C^1[0, \tau]$ into their complete spline interpolants $\mathcal{C}[g]$ (based on the fixed knot sequence $0 = a_0 < a_1 \ldots < a_k = \tau$); thus in this section our $\mathcal{C}$ is de Boor's operator $I_4$. Thus we have

$$\mathbb{H}_{n,k} = \mathcal{C}[\mathbb{Y}_n], \qquad H_{n,k} = \mathcal{C}[Y].$$

It follows from the formula for $c_{4,i}$ in (5) on page 40 of de Boor [5] that the slope of $\mathbb{H}_{n,k}^{(2)}$ on the interval $[a_{j-1}, a_j]$ is given by

$$B_j \equiv B_j(CS) = \frac{12}{(\Delta_j a)^3} \left( \frac{\mathbb{H}_{n,k}^{(1)}(a_{j-1}) + \mathbb{H}_{n,k}^{(1)}(a_j)}{2} \Delta_j a - \Delta_j \mathbb{Y}_n \right)$$

where $\Delta_j f \equiv f(a_j) - f(a_{j-1})$ for $j = 1, \ldots, k$ and any function $f$ on $[0, \tau]$.

In the following we will let $\mathcal{H}$ denote the Hermite interpolation operator that maps $\mathbb{Y}_n$ to $\mathbb{H}_n$: thus $\mathbb{H}_{n,Herm} = \mathcal{H}[\mathbb{Y}_n]$, $\mathbb{H}_{n,Herm}^{(1)} = (\mathcal{H}[\mathbb{Y}_n])^{(1)}$, and so forth. It is important to note that the corresponding slopes of the second derivative of the Hermite interpolant, $\mathbb{H}_{n,Herm}^{(2)} = (\mathcal{H}[\mathbb{Y}_n])^{(2)}$ on $[a_{j-1}, a_j]$ are given by the same formula as in the last display, but with $\mathbb{H}_{n,k}^{(1)}(a_i)$ replaced by $\mathbb{Y}_n^{(1)}(a_i) = \mathbb{F}_n(a_i)$, $i = j-1, j$:

$$(13) \qquad \widetilde{B}_j \equiv B_j(Herm) = \frac{12}{(\Delta_j a)^3} \left( \frac{\mathbb{F}_n(a_{j-1}) + \mathbb{F}_n(a_j)}{2} \Delta_j a - \Delta_j \mathbb{Y}_n \right).$$

Note that $\widetilde{B}_j$ is expressed explicitly as a function of the data via $\mathbb{F}_n$ and $\mathbb{Y}_n$, whereas $B_j$ still involves $\mathbb{H}_{n,k} = \mathcal{C}[\mathbb{Y}_n]$ and hence also the interpolation operator $\mathcal{C}$. Ordering of the slopes $\widetilde{B}_j$ can be shown using only Lemma 3.1 and Lemma 4.5, but (unfortunately) the generalization of Marshall's lemma does not apply to the Hermite interpolant because the second derivative $\mathbb{H}_{n,Herm}^{(2)}$ is not continuous at the knots. This last formula (13) agrees with the formulas for $H$ and $H_n$ in Groeneboom, Jongbloed and Wellner [10] and Groeneboom, Jongbloed and Wellner [11]; in particular (13) can be viewed as a finite sample analogue of the 3rd derivative of the interpolant $H$ given in Groeneboom, Jongbloed and Wellner [10], page 1631, but based on the *fixed knots* $\{a_j\}$ rather than random knots determined by the optimization procedure. Note that the least squares estimator $\widetilde{f}_n = \widetilde{\mathbb{H}}_n^{(2)}$ can be viewed as the second derivative of either the Hermite interpolant or the complete cubic spline interpolant of $\mathbb{Y}_n$ since these two interpolants have been forced equal by the optimization procedure which determines the knots as random functions of the data.

Set

$$A_n \equiv \left\{ \mathbb{H}_{n,k_n}^{(2)} \text{ is convex on } [0,\tau] \right\} = \bigcap_{j=1}^{k-1} \{B_j \le B_{j+1}\}.$$

To prove **B**, we want to bound

$$P(A_n^c) \le \sum_{j=1}^{k-1} P(B_j > B_{j+1}).$$

To prepare for this, we define

$$T_{n,j} = \frac{(\mathcal{C}[\mathbb{Y}_n])^{(1)}(a_{j-1}) + (\mathcal{C}[\mathbb{Y}_n])^{(1)}(a_j)}{2}\Delta_j a - \Delta_j \mathbb{Y}_n,$$

$$R_{n,j} = \frac{\mathbb{Y}_n^{(1)}(a_{j-1}) + \mathbb{Y}_n^{(1)}(a_j)}{2}\Delta_j a - \Delta_j \mathbb{Y}_n,$$

$$t_{n,j} = \frac{(\mathcal{C}[Y])^{(1)}(a_{j-1}) + (\mathcal{C}[Y])^{(1)}(a_j)}{2}\Delta_j a - \Delta_j Y,$$

$$r_{n,j} = \frac{Y^{(1)}(a_{j-1}) + Y^{(1)}(a_j)}{2}\Delta_j a - \Delta_j Y.$$

We will frequently suppress the dependence of all of these quantities on $n$, and simply write $T_j$ for $T_{n,j}$, $R_j$ for $R_{n,j}$, and so forth. Now $B_j = 12T_j/(\Delta_j a)^3$, $\tilde{B}_j = 12R_j/(\Delta_j a)^3$, and we can write

(14) $$T_j - r_j = T_j - t_j + t_j - r_j$$
$$= R_j - r_j + \{T_j - t_j - (R_j - r_j)\} + t_j - r_j$$
(15) $$\equiv R_j - r_j + W_j + b_j.$$

We regard $R_j - r_j$ as the main random term to be controlled, and view $T_j - t_j - (R_j - r_j) \equiv W_j$ and $t_j - r_j \equiv b_j$ as second order terms, the last of which is deterministic. Thus our strategy will be to first develop an appropriate exponential bound for $|R_j - r_j|$, and then by further separate bounds for $W_j$ and $b_j$, derive an exponential bound for $|T_j - r_j|$.

For $0 \leq s < t < \infty$, define the family of functions $h_{s,t}$ by

$$h_{s,t}(x) = (x - (s+t)/2)1_{(s,t]}(x).$$

Note that

$$Ph_{s,t} = \frac{1}{2}(F(t) + F(s))(t - s) - \int_s^t F(u)du,$$

$$\mathbb{P}_n h_{s,t} = \frac{1}{2}(\mathbb{F}_n(t) + \mathbb{F}_n(s))(t - s) - \int_s^t \mathbb{F}_n(u)du,$$

and, furthermore,

$$r_j = Ph_{a_{j-1},a_j}, \qquad R_j = \mathbb{P}_n h_{a_{j-1},a_j}.$$

Here is a (partial) analogue of Kiefer and Wolfowitz's Lemma 1.

**Lemma 3.1.** *Suppose that $\tilde{\gamma}_1(F) < \infty$ and $R < \infty$. Let $h_{s,t}(x) = (x - (s+t)/2)1_{(s,t]}(x)$, $s = a_{j-1}^{(k)} \equiv a_{j-1}$, and $t = a_j^{(k)} \equiv a_j$ so that $t - s = a_j - a_{j-1} = k^{-1}(1/f(a_j^*))$ for some $a_j^* \in [a_{j-1}, a_j]$. Then if $\delta_n \to 0$ and $k \geq 5\tilde{\gamma}_1(F)R$,*

$$Pr(|R_j - r_j| > \delta_n p_n^3) = Pr(|\mathbb{P}_n - P|(h_{s,t}) > \delta_n p_n^3)$$
$$= 2\exp\left(-\frac{3n\delta_n^2 f^2(a_j^*)p_n^3}{1 + p_n \delta_n f(a_j^*)}\right)$$
$$\leq 2\exp(-3n\delta_n^2 p_n^3 f^2(a_j^*)(1 + o(1)))$$

*where $o(1)$ depends on $f(a_j^*)$, $k_n$, and $\delta_n$.*

*Proof.* First note that $|h_{s,t}|$ is bounded by $(t-s)/2$. Thus by Bernstein's inequality (see e.g. van der Vaart and Wellner [23], page 102),

$$Pr(|\mathbb{P}_n h_{s,t} - P h_{s,t}| > x) \leq 2 \exp\left(-\frac{nx^2/2}{\sigma^2 + Mx/3}\right)$$

for $\sigma^2 \geq Var_F(h_{s,t}(X))$, $M = (t-s)/2 = 1/(2f(a_j^*)k) = [1/(2f(a_j^*))]p_n$, and $x > 0$. Note that

$$\begin{aligned}
Var(h_{s,t}(X)) \leq E h_{s,t}^2(X) &= \int_s^t (x - (t+s)/2)^2 dF(x) \\
&\leq f(s)(t-s)^3/12 \\
&= f(s)k^{-3}/(12f^3(a_j^*)) = f(s)p_n^3/(12f^3(a_j^*)) \\
&\leq p_n^3/(6f^2(a_j^*))
\end{aligned}$$

for $k \geq 5\tilde{\gamma}_1(F)R$ by Lemma 4.1. Then we obtain

$$\begin{aligned}
Pr(|\mathbb{P}_n h_{s,t} - P h_{s,t}| &> \delta_n p_n^3) \\
&\leq 2 \exp\left(-\frac{n\delta_n^2 p_n^6/2}{p_n^3/(6f(a_j^*)^2) + p_n \delta_n p_n^3/(6f(a_j^*))}\right) \\
&= 2 \exp\left(-\frac{n\delta_n^2 f^2(a_j^*) p_n^3}{1/3 + p_n f(a_j^*)\delta_n/3}\right) \\
&= 2 \exp\left(-\frac{3n\delta_n^2 f^2(a_j^*) p_n^3}{1 + p_n \delta_n f(a_j^*)}\right) \\
&= 2 \exp\left(-3n\delta_n^2 p_n^3 f^2(a_j^*)(1 + o(1))\right)
\end{aligned}$$

where the $o(1)$ term depends on $f(t) = f(a_{j+1})$, $p_n = 1/k_n$, and $\delta_n$.   □

**Remark.** Note that taking $\delta_n = C/k_n$ in Lemma 3.1 yields

$$Pr(|\mathbb{P}_n h_{s,t} - P h_{s,t}| > C p_n^4) \leq 2 \exp(-3(nC^2 f^2(a_j^*)/k_n^5)(1 + o(1)))$$

which seems quite analogous to Lemma 4 of Kiefer and Wolfowitz (1976), but with the power of 3 replaced by 5.

The following lemma gives a more complete version of Lemma 3.1 in that it provides an exponential bound for $|T_j - r_j|$.

**Lemma 3.2.** *Suppose that the hypotheses of Theorem 2.2 hold: $\beta_2(F, \tau) < \infty$, $\gamma_2(F, \tau) < \infty$, $\tilde{\gamma}_1(F) < \infty$ and $R < \infty$. Then if $\delta_n = Cp_n$ for some constant $C$ and $k \geq \{5R \vee 3\}\tilde{\gamma}_1(F)$,*

$$Pr(|T_j - r_j| > 3\delta_n p_n^3) \leq 6 \exp\left(-\frac{(100)^{-1} n \delta_n^2 f^2(a_j^*) p_n^3}{1 + 30^{-1} p_n \delta_n f(a_j^*)}\right).$$

*Proof.* This follows from a combination of Lemma 3.1, Lemma 4.2, and Lemma 4.3. Lemma 4.2 yields

$$|b_j| \equiv |t_j - r_j| \leq R^4 o(1) p_n^4 \leq \delta_n p_n^3$$

if $n$ (and hence $k_n$) is sufficiently large. This implies that

$$\begin{aligned}
Pr(|T_j - r_j| > 3\delta_n p_n^3) &\leq Pr(|T_j - t_j| > 3\delta_n p_n^3 - |t_j - r_j|) \\
&\leq Pr(|T_j - t_j| > 2\delta_n p_n^3).
\end{aligned}$$

In view of the decomposition (15), this yields

$$Pr(|T_j - r_j| > 3\delta_n p_n^3) \leq Pr(|R_j - r_j| > \delta_n p_n^3) + Pr(|A_j| > \delta_n p_n^3)$$

$$\leq 6 \exp\left(-\frac{(100)^{-1} n \delta_n^2 f^2(a_j^*) p_n^3}{1 + 30^{-1} p_n \delta_n f(a_j^*)}\right)$$

by Lemma 3.1, Lemma 4.3, and the fact that

$$\frac{100^{-1} A}{1 + 30^{-1} B} = \frac{30}{100} \frac{A}{30 + B} \leq \frac{3A}{1 + B}$$

for $A, B > 0$. $\square$

**Lemma 3.3.** *Suppose that* $\beta_2 \equiv \beta_2(F, \tau) > 0$, $\tilde{\gamma}_1 \equiv \tilde{\gamma}_1(F, \tau) < \infty$ *and* $R \equiv R(f, \tau) < \infty$ *for some* $\tau < \alpha_1(F) \equiv \inf\{t : F(t) = 1\}$. *Let*

$$A_n \equiv \{\mathbb{H}_{n,k_n}^{(2)} \text{ is convex on } [0, \tau]\}.$$

*Then*

(16) $$P(A_n^c) \leq 12 k_n \exp\left(-K \beta_2^2(F, \tau) n p_n^5\right)$$

*where* $K^{-1} = 8^2 \cdot 144^2 \cdot 16 \cdot 200 = 4{,}246{,}732{,}800 \leq 4.3 \cdot 10^9$.

*Proof.* Since

$$A_n^c \equiv \bigcup_{j=1}^{k_n - 1} \{B_j > B_{j+1}\},$$

it follows that

$$P(A_n^c) \leq \sum_{j=1}^{k_n-1} P(B_j > B_{j+1})$$

$$= \sum_{j=1}^{k_n-1} P\left(B_j > B_{j+1}, |T_i - r_i| \leq 3\delta_{n,j} p_n^3, \ i = j, j+1\right)$$

$$+ \sum_{j=1}^{k_n-1} P\left(B_j > B_{j+1}, |T_i - r_i| > 3\delta_{n,j} p_n^3 \text{ for } i = j \text{ or } i = j+1\right)$$

$$\leq \sum_{j=0}^{m_n-1} P\left(B_j > B_{j+1}, |T_i - r_i| \leq 3\delta_{n,j} p_n^3, \ i = j, j+1\right)$$

$$+ \sum_{j=0}^{k_n-1} \left\{P\left(|T_j - r_j| > 3\delta_{n,j} p_n^3\right) + P\left(|T_{j+1} - r_{j+1}| > 3\delta_{n,j} p_n^3\right)\right\}$$

(17) $$= I_n + II_n$$

where we take

$$\delta_{n,j} = \frac{C(F, \tau)}{k_n f(a_j^*)} = p_n \frac{C(F, \tau)}{f(a_j^*)} \equiv \frac{\delta_n}{f(a_j^*)};$$

here $a_j^* \in [a_{j-1}, a_j]$ satisfies $\Delta_j a = a_j - a_{j-1} = 1/(k_n f(a_j^*))$, and $C(F, \tau)$ is a constant to be determined. We first bound $II_n$ from above. By Lemma 3.2, we know that

$$P\left(|T_j - r_j| > 3\delta_{n,j} p_n^3\right) \leq 6 \exp\left(-\frac{(100)^{-1} n \delta_{n,j}^2 f^2(a_j^*) p_n^3}{1 + 30^{-1} p_n \delta_{n,j} f(a_j^*)}\right)$$

where $\delta_{n,j}^2 f^2(a_j^*) p_n^3 = C^2(F, \tau) p_n^5$ and

$$\frac{1}{1 + 30^{-1} p_n \delta_{n,j} f(a_j^*)} = \frac{1}{1 + 30^{-1} C(F, \tau) p_n^2} > \frac{1}{2}$$

when $k_n > [30^{-1} C(F, \tau)]^{1/2}$. Hence,

(18) $$P\left(|T_j - r_j| > 3\delta_{n,j} p_n^3\right) \leq 6 \exp\left(-200^{-1} C^2(F, \tau) n p_n^5\right).$$

We also have

$$P\left(|T_{j+1} - r_{j+1}| > 3\delta_{n,j} p_n^3\right) \leq 6 \exp\left(-\frac{100^{-1} n \delta_{n,j}^2 f^2(a_{j+1}^*) p_n^3}{1 + 30^{-1} p_n \delta_{n,j} f(a_{j+1}^*)}\right)$$

where $a_{j+1}^* \in [a_j, a_{j+1}]$ and $a_{j+1} - a_j = \Delta_{j+1} a = 1/(k_n f(a_{j+1}^*))$. By Lemma 5.1 we have $f(a_j)/f(a_{j+1}) \leq 2$ if $k_n \geq 5\tilde{\gamma}_1(F, \tau)R$. But this implies that $f(a_j^*)/f(a_{j+1}^*) \leq 4$ since

$$\frac{f(a_j^*)}{f(a_{j+1}^*)} = \frac{f(a_j^*)}{f(a_j)} \cdot \frac{f(a_j)}{f(a_{j+1})} \cdot \frac{f(a_{j+1})}{f(a_{j+1}^*)}$$
$$\leq \frac{f(a_{j-1})}{f(a_j)} \frac{f(a_j)}{f(a_{j+1})} \frac{f(a_{j+1})}{f(a_{j+1}^*)} \leq 2 \cdot 2 \cdot 1 = 4.$$

Hence, we can write

$$\delta_{n,j}^2 f^2(a_{j+1}^*) = \frac{1}{k_n^2} C^2(F, \tau) \frac{f^2(a_{j+1}^*)}{f^2(a_j^*)} \geq \frac{1}{k_n^2} C^2(F, \tau) \frac{1}{16} = \frac{C^2(F, \tau)}{16} p_n^2$$

and, since $f(a_{j+1}^*)/f(a_j^*) \leq 1$,

$$\frac{1}{1 + 30^{-1} p_n \delta_{n,j} f(a_{j+1}^*)} = \frac{1}{1 + 30^{-1} C(F, \tau) p_n^2 f(a_{j+1}^*)/f(a_j^*)}$$
$$\geq \frac{1}{1 + 30^{-1} C(F, \tau) p_n^2} > \frac{1}{2}$$

when $k_n > [30^{-1} C(F, \tau)]^{1/2}$. Thus, we conclude that

(19) $$P\left(|T_{j+1} - r_{j+1}| > 3\delta_{n,j} p_n^3\right) \leq 6 \exp\left(-\frac{200^{-1}}{16} C^2(F, \tau) n p_n^5\right)$$

Combining (18) and (19), we get

$$II_n \leq 12 k_n \exp\left(-\frac{200^{-1}}{16} C^2(F, \tau) n p_n^5\right).$$

Now we need to handle $I_n$. Recall that

$$B_j = 12\frac{T_j}{(\Delta_j a)^3}, \qquad B_{j+1} = 12\frac{T_{j+1}}{(\Delta_{j+1} a)^3}.$$

Thus, the event

$$\left\{B_j > B_{j+1}, \ |T_i - r_i| \leq 3\delta_{n,j} p_n^3, \ i = j, j+1\right\}$$

is equal to the event

$$\left\{\frac{T_j}{(\Delta_j a)^3} > \frac{T_{j+1}}{(\Delta_{j+1} a)^3}, \ |T_i - r_i| \leq 3\delta_{n,j} p_n^3, \ i = j, j+1\right\}.$$

Then, it follows that

$$\frac{T_j}{(\Delta_j a)^3} \leq \frac{r_j}{(\Delta_j a)^3} + \frac{3\delta_{n,j} p_n^3}{(\Delta_j a)^3}$$

and

$$\frac{T_{j+1}}{(\Delta_{j+1} a)^3} \geq \frac{r_{j+1}}{(\Delta_{j+1} a)^3} - \frac{3\delta_{n,j} p_n^3}{(\Delta_{j+1} a)^3},$$

and hence

$$\begin{aligned}\frac{T_j}{(\Delta_j a)^3} &\leq \left[\frac{r_j}{(\Delta_j a)^3} - \frac{r_{j+1}}{(\Delta_{j+1} a)^3}\right] + \left[\frac{r_{j+1}}{(\Delta_{j+1} a)^3} - \frac{3\delta_{n,j} p_n^3}{(\Delta_{j+1} a)^3}\right] \\ &\quad + \left[\frac{3\delta_{n,j} p_n^3}{(\Delta_j a)^3} + \frac{3\delta_{n,j} p_n^3}{(\Delta_{j+1} a)^3}\right] \\ &\leq \left[\frac{r_j}{(\Delta_j a)^3} - \frac{r_{j+1}}{(\Delta_{j+1} a)^3}\right] + \frac{T_{j+1}}{(\Delta_{j+1} a)^3} \\ &\quad + \left[\frac{3\delta_{n,j} p_n^3}{(\Delta_j a)^3} + \frac{3\delta_{n,j} p_n^3}{(\Delta_{j+1} a)^3}\right].\end{aligned}$$

The first term in the right side of the previous inequality is the leading term in the sense that it determines the sign of the difference of the slope of $\mathbb{H}_{n,k_n}^{(2)}$. By Lemma 4.5, we can write

$$\frac{r_j}{(\Delta_j a)^3} - \frac{r_{j+1}}{(\Delta_{j+1} a)^3} \leq -\frac{1}{12} f''(a_j^{**})\Delta_j a + \frac{1}{24}(\overline{f}_j''\Delta_j a - \underline{f}_{j+1}''\Delta_{j+1} a).$$

Let $a_j^* \in [a_{j-1}, a_j]$ such that $\Delta_j a = p_n[f(a_j^*)]^{-1}$. Then, we can write

$$\begin{aligned}&\frac{3\delta_{n,j} p_n^3}{(\Delta_j a)^3} + \frac{3\delta_{n,j} p_n^3}{(\Delta_{j+1} a)^3} - \frac{1}{12}f''(a_j^{**})\Delta_j a + \frac{1}{24}(\overline{f}_j''\Delta_j a - \underline{f}_{j+1}''\Delta_{j+1} a) \\ &\leq 6\delta_{n,j} f^3(a_j^*) - \frac{1}{12} f''(a_j^{**})\Delta_j a + \frac{1}{24}(\overline{f}_j''\Delta_j a - \underline{f}_{j+1}''\Delta_{j+1} a) \\ &= 6f^2(a_j^*)\left\{\delta_n - \frac{1}{72}\frac{f''(a_j^{**})}{f^3(a_j^*)} p_n + \frac{1}{144 f^2(a_j^*)}(\overline{f}_j''\Delta_j a - \underline{f}_{j+1}''\Delta_{j+1} a)\right\} \\ &= 6f^2(a_j^*)\left\{\delta_n - \frac{1}{72}\frac{f''(a_j^{**})}{f^3(a_j^*)} p_n + \frac{1}{144 f^3(a_j^*)}\left(\overline{f}_j'' - \underline{f}_{j+1}''\frac{\Delta_{j+1} a}{\Delta_j a}\right) p_n\right\} \\ &= 6f^2(a_j^*)\left\{\delta_n - \frac{1}{72}\frac{f''(a_j^{**})}{f^3(a_j^*)} p_n + \frac{1}{144}\frac{\underline{f}_{j+1}''}{f^3(a_j^*)}\left(\frac{\overline{f}_j''}{\underline{f}_{j+1}''} - \frac{\Delta_{j+1} a}{\Delta_j a}\right) p_n\right\}\end{aligned}$$

$$= 6f^2(a_j^*)\left\{\delta_n - \frac{1}{72}\frac{f''(a_j^{**})}{f^3(a_j^{**})}\frac{f^3(a_j^{**})}{f^3(a_j^*)}p_n\right.$$
$$\left. + \frac{1}{144}\frac{\underline{f}''_{j+1}}{f^3(a_j^*)}\left(\frac{\overline{f}''_j}{\underline{f}''_{j+1}} - \frac{\Delta_{j+1}a}{\Delta_j a}\right)p_n\right\}$$
$$\leq 6f^2(a_j^*)\left\{\delta_n - \frac{1}{72}\frac{\beta_2(F,\tau)}{8}p_n + \frac{1}{144}\frac{\underline{f}''_{j+1}}{f^3(a_j^*)}\left(\frac{\overline{f}''_j}{\underline{f}''_{j+1}} - \frac{\Delta_{j+1}a}{\Delta_j a}\right)p_n\right\}$$
$$= 6f^2(a_j^*)\left\{\delta_n - \frac{1}{72}\frac{\beta_2(F,\tau)}{8}p_n\right.$$
$$\left. + \frac{1}{144}\frac{\underline{f}''_{j+1}}{f^3(a_j^*)}\left(\frac{\overline{f}''_j}{\underline{f}''_{j+1}} - 1 + 1 - \frac{\Delta_{j+1}a}{\Delta_j a}\right)p_n\right\}$$

where (using arguments similar to those of Lemma 4.2 and taking the bound on $|\overline{f}''_j - \underline{f}''_{j+1}|$ to be $\epsilon\|f''\|$ which is possible by uniform continuity of $f''$ on $[0,\tau]$)

$$\frac{\overline{f}''_j}{\underline{f}''_{j+1}} - 1 \leq \left|\frac{\overline{f}''_j}{\underline{f}''_{j+1}} - 1\right| \leq \frac{\epsilon f(\tau)^3\gamma_2(F,\tau)}{\underline{f}''_{j+1}}$$

if $k_n > \max(5\tilde{\gamma}_1(F,\tau)R, (\sqrt{2}+1)R/\eta)$ for a given $\eta > 0$ and

$$1 - \frac{\Delta_{j+1}a}{\Delta_j a} \leq \left|\frac{\Delta_{j+1}a}{\Delta_j a} - 1\right| \leq 8\tilde{\gamma}_1(F,\tau)p_n.$$

Hence

$$\frac{3\delta_{n,j}p_n^3}{(\Delta_j a)^3} + \frac{3\delta_{n,j}p_n^3}{(\Delta_{j+1}a)^3} - \frac{1}{12}f''(a_j^{**})\Delta_j a + \frac{1}{24}(\overline{f}''_j\Delta_j a - \underline{f}''_{j+1}\Delta_{j+1}a)$$
$$\leq 6f^2(a_j^*)\left\{\delta_n - \frac{1}{72}\frac{\beta_2(F,\tau)}{8}p_n + \frac{1}{144}\epsilon\gamma_2(F,\tau)\, p_n\right.$$
$$\left. + \frac{8}{144}\gamma_2(F,\tau)\tilde{\gamma}_1(F,\tau)p_n^2\right\}$$

where we can choose $\epsilon$ and $p_n$ small enough so that

$$\frac{1}{144}\epsilon\gamma_2(F,\tau) + \frac{8}{144}\gamma_2(F,\tau)\tilde{\gamma}_1\, p_n \leq \frac{1}{2\cdot 72\cdot 8}\beta_2(F,\tau);$$

for example

$$\epsilon < \frac{1}{16}\frac{\beta_2(F,\tau)}{\gamma_2(F,\tau)}, \quad k_n = p_n^{-1} > 16\cdot 8\frac{\tilde{\gamma}_1(F,\tau)}{\beta_2(F,\tau)}.$$

The above choice yields

$$\frac{3\delta_{n,j}p_n^3}{(\Delta_j a)^3} + \frac{3\delta_{n,j}p_n^3}{(\Delta_{j+1}a)^3} - \frac{1}{12}\underline{f}''_j\Delta_j a + \frac{1}{24}(\overline{f}''_j\Delta_j a - \underline{f}''_{j+1}\Delta_{j+1}a)$$
$$\leq 6f^2(a_j^*)\left\{\delta_n - \frac{\beta_2(F,\tau)}{8\cdot 144}p_n\right\} = 0$$

by choosing
$$\delta_n = C(F,\tau)p_n = \frac{\beta_2(F,\tau)}{8 \cdot 144}p_n;$$
i.e. $C(F,\tau) = \beta_2(F,\tau)/(8 \cdot 144)$. For such a choice, the first term $I_n$ in (17) is identically equal to 0. □

## 4. Appendix 1: technical lemmas

**Lemma 4.1.** *Under the hypotheses of Theorem 2.2,*
$$1 \le \frac{f(a_{j-1})}{f(a_j)} \le \frac{\Delta_{j+1}a}{\Delta_j a} \le 2$$
*uniformly in $j$ if $k \ge 5\tilde{\gamma}_1 R$.*

*Proof.* Note that for each interval $I_j = [a_{j-1}, a_j]$ we have
$$p_n = \int_{I_j} f(x)dx = f(a_j^*)\Delta_j a \begin{cases} \ge f(a_j)\Delta_j a \\ \le f(a_{j-1})\Delta_j a \end{cases}$$
where $a_j^* \in I_j$. Thus
$$\frac{p_n}{\Delta_{j+1}a} \le f(a_j) \le \frac{p_n}{\Delta_j a}$$
and
$$\frac{p_n}{\Delta_j a} \le f(a_{j-1}) \le \frac{p_n}{\Delta_{j-1}a}.$$
It follows that
$$1 \le \frac{f(a_{j-1})}{f(a_j)} \le \frac{\Delta_{j+1}a}{\Delta_{j-1}a} = \frac{\Delta_{j+1}a}{\Delta_j a}\frac{\Delta_j a}{\Delta_{j-1}a}.$$
Thus we will establish a bound for $\Delta_{j+1}a/\Delta_j a$. Note that with $c \equiv F(\tau) < 1$
$$\begin{aligned}
\Delta_{j+1}a = a_{j+1} - a_j &= F^{-1}(\frac{j+1}{k}c) - F^{-1}(\frac{j}{k}c) \\
&= \frac{c}{k}\frac{1}{f(a_j)} + \frac{c^2}{2k^2}\frac{-f'(\xi_{j+1})}{f^3(\xi_{j+1})} \\
&= \frac{c}{k}\frac{1}{f(a_j)}\left\{1 + \frac{c}{2k}\frac{-f'(\xi_{j+1})}{f^2(\xi_{j+1})}\frac{f(a_j)}{f(\xi_{j+1})}\right\} \\
&\le \frac{c}{k}\frac{1}{f(a_j)}\left\{1 + \frac{c\tilde{\gamma}_1}{2k}R\right\}.
\end{aligned}$$
for some $\xi_{j+1} \in I_{j+1}$, where $\xi_{j+1} \in I_{j+1}$, $R < \infty$, and $\tilde{\gamma}_1 < \infty$.

Similarly, expanding to second order (about $a_j$ again!),
$$\begin{aligned}
\Delta_j a = a_j - a_{j-1} &= F^{-1}(\frac{j}{k}c) - F^{-1}(\frac{j-1}{k}c) \\
&= \frac{c}{k}\frac{1}{f(a_j)} + \frac{c^2}{2k^2}\frac{f'(\xi_j)}{f^3(\xi_j)} \\
&= \frac{c}{k}\frac{1}{f(a_j)}\left\{1 + \frac{c}{2k}\frac{f'(\xi_j)}{f^2(\xi_j)}\frac{f(a_j)}{f(\xi_j)}\right\}
\end{aligned}$$

$$\geq \frac{c}{k}\frac{1}{f(a_j)}\left\{1 + \frac{c}{2k}\frac{f'(\xi_j)}{f^2(\xi_j)}\right\}$$

since $f(a_j)/f(\xi_j) \leq 1$ and $f'(\xi_j) < 0$

$$\geq \frac{c}{k}\frac{1}{f(a_j)}\left\{1 - \frac{c\tilde{\gamma}_1}{2k}\right\}.$$

where $\xi_j \in I_j$. Thus it follows that for $k = k_n$ so large that $\tilde{\gamma}_1/(2k) \leq 1/2$ we have

$$\frac{\Delta_{j+1}a}{\Delta_j a} \leq \frac{1 + \frac{c\tilde{\gamma}_1}{2k}R}{1 - \frac{c\tilde{\gamma}_1}{2k}}$$

$$\leq \left(1 + \frac{c\tilde{\gamma}_1}{2k}R\right)\left(1 + \frac{c\tilde{\gamma}_1}{k}\right)$$

$$= 1 + \frac{c\tilde{\gamma}_1}{k}(R/2 + 1) + \frac{c^2\tilde{\gamma}_1^2}{2k^2}R$$

$$< 1 + \frac{\tilde{\gamma}_1(R+1)}{k}$$

if $k = k_n \geq \tilde{\gamma}_1$. The last inequality here follows from

$$\frac{\tilde{\gamma}_1}{k}(R/2 + 1) + \frac{\tilde{\gamma}_1^2}{2k^2}R \leq \frac{\tilde{\gamma}_1}{k}(R + \alpha)$$

if and only if

$$(R/2 + 1) + \frac{\tilde{\gamma}_1}{2k}R \leq R + \alpha$$

or, equivalently, if and only if

$$\frac{\tilde{\gamma}_1}{2k}R \leq R/2 + \alpha - 1, \quad \text{or} \quad k \geq \tilde{\gamma}_1\frac{R}{R + 2(\alpha - 1)} = \tilde{\gamma}_1$$

if $\alpha = 1$. It now follows that

$$1 \leq \frac{f(a_{j-1})}{f(a_j)} \leq \frac{\Delta_{j+1}a}{\Delta_{j-1}a} = \frac{\Delta_{j+1}a}{\Delta_j a}\frac{\Delta_j a}{\Delta_{j-1}a\Delta_{j-1}a} \leq 2$$

if

$$\frac{\Delta_{i+1}a}{\Delta_i a} \leq \sqrt{2}$$

for $i = j-1, j$. But these inequalities hold if $k$ is so large that $1 + \frac{\tilde{\gamma}_1(R+1)}{k} \leq \sqrt{2}$, or $k \geq 5\tilde{\gamma}_1 R \geq \tilde{\gamma}_1(R+1)/(\sqrt{2}-1)$ since $R \geq 1$ and $1/(\sqrt{2}-1) \leq 5/2$. $\square$

**Lemma 4.2.** *Under the hypotheses of Theorem 2.2,*

$$\frac{|t_j - r_j|}{(\Delta_j a)^4} = o(1)$$

*where the $o(1)$ depends only on $\tau$, $\tilde{\gamma}_1(F,\tau)$, and $\gamma_2(F,\tau)$.*

**Remark.** Note that

(20) $$\max_{1 \leq j \leq k}|t_j - r_j| \leq \frac{1}{24}|a|^4\|Y^{(4)}\| = \frac{1}{24}|a|^4\|f''\| \leq \frac{1}{24}R\gamma_2(F)p_n^4.$$

This follows since

$$r_j - t_j = \frac{1}{2}\left(Y^{(1)}(a_{j-1}) + Y^{(1)}(a_j)\right.$$
$$\left. - ((\mathcal{C}[Y])^{(1)}(a_{j-1}) - (\mathcal{C}[Y])^{(1)}(a_j))\right)\Delta_j a$$
$$= \frac{1}{2}\left\{\left(Y^{(1)}(a_{j-1}) - (\mathcal{C}[Y])^{(1)}(a_{j-1})\right)\right.$$
$$\left. + \left(Y^{(1)}(a_j) - (\mathcal{C}[Y])^{(1)}(a_j)\right)\right\}\Delta_j a,$$

and hence from de Boor [5], (20), page 56, it follows that

$$|r_j - t_j| \leq \frac{1}{24}|a|^3 \|Y^{(4)}\|\Delta_j a \leq \frac{1}{24}|a|^4 \|f^{(2)}\|,$$

and this yields (20). The claim of Lemma 4.2 is stronger because it makes a statement about the differences $t_j - r_j$ relative to $(\Delta_j a)^4$; this is possible because only differences between the derivative of the derivative of $Y$ and the derivative of its interpolant $\mathcal{C}[Y]$ *at the knots* are involved.

*Proof.* We have

$$(21) \qquad r_j - t_j = \frac{1}{2}\bigl(\mathcal{E}_Y^{(1)}(a_{j-1}) + \mathcal{E}_Y^{(1)}(a_j)\bigr)\Delta_j a,$$

where $\mathcal{E}_g = g - \mathcal{C}[g]$. Now, using the result of Problem 2a, Chapter V of de Boor [5] (compare also with the formula (3.52) given in Nürnberger [20]), we have

$$\delta_j \mathcal{E}_Y^{(1)}(a_{j-1}) + 2\mathcal{E}_Y^{(1)}(a_j) + (1-\delta_j)\mathcal{E}_Y^{(1)}(a_{j+1}) = \beta_j$$

for $j = 0, \cdots, k-1$, where

$$\delta_j = \frac{a_{j+1} - a_j}{a_{j+1} - a_{j-1}} = \frac{\Delta_{j+1} a}{\Delta_j a + \Delta_{j+1} a}$$

and

$$\beta_j = \frac{\delta_j(-\Delta_j a)^3 f''(\xi_{1,j}) + (1-\delta_j)(\Delta_{j+1} a)^3 f''(\xi_{2,j})}{24},$$

$\xi_{1,j}, \xi_{2,j} \in [a_{j-1}, a_{j+1}]$. By Problem IV 7(a) in de Boor [5] and the techniques used in Chapter III (see in particular equation (9)), a bound on the maximal value at the knots of the derivative interpolation error can be derived using the following inequality

$$(22) \qquad \max_{0 \leq j \leq k} |\mathcal{E}_Y^{(1)}(a_j)| \leq \max\left(|\mathcal{E}_Y^{(1)}(a_0)|, \max_{1 \leq j \leq k-1} |\beta_j|, |\mathcal{E}_Y^{(1)}(a_k)|\right).$$

By definition of the complete cubic spline, $\mathcal{E}_Y^{(1)}(a_0) = \mathcal{E}_Y^{(1)}(a_k) = 0$. Thus, we will focus now on getting a sharp bound for $\max_{1 \leq j \leq k-1} |\beta_j|$ under our hypotheses. This will be achieved as follows:

- Expanding $\delta_j$ around $1/2$: We have

$$\delta_j = \frac{a_{j+1} - a_j}{(a_{j+1} - a_j) + (a_j - a_{j-1})} = \frac{k_n^{-1}[f(a_{j+1}^*)]^{-1}}{k_n^{-1}[f(a_{j+1}^*)]^{-1} + k_n^{-1}[f(a_j^*)]^{-1}},$$

where $a_j^* \in [a_{j-1}, a_j]$ and $a_j^{**} \in [a_j, a_{j+1}]$, and hence

$$\delta_j = \frac{1}{2} + \frac{f(a_{j+1}^*) - f(a_j^*)}{2(f(a_j^*) + f(a_{j+1}^*))}$$

$$= \frac{1}{2} + \frac{f'(a_j^{**})}{2(f(a_j^*) + f(a_{j+1}^*))}(a_j^{**} - a_j^*)$$

$$= \frac{1}{2} + \frac{f'(a_j^{**})}{2(f(a_j^*) + f(a_{j+1}^*))} \frac{a_j^{**} - a_j^*}{a_j - a_{j-1}} \Delta_j a = \frac{1}{2} + M_j \, \Delta_j a$$

where

$$|M_j| = \left| \frac{f'(a_j^{**})}{2(f(a_j^*) + f(a_{j+1}^*))} \frac{a_j^{**} - a_j^*}{a_j - a_{j-1}} \right|$$

$$\leq \frac{|f'(a_{j-1})|}{4f(a_{j+1})} \frac{a_{j+1} - a_{j-1}}{a_j - a_{j-1}}$$

$$\leq \frac{|f'(a_{j-1})|}{4f(a_{j-1})} \frac{f(a_{j-1})}{f(a_{j+1})} \left( \frac{a_{j+1} - a_j}{a_j - a_{j-1}} + 1 \right)$$

$$\leq \frac{|f'(a_{j-1})|}{4f(a_{j-1})} 2 \cdot 2 \cdot (\sqrt{2} + 1), \quad \text{for } k_n > 5\tilde{\gamma}_1 R$$

$$= (\sqrt{2} + 1) \frac{|f'(a_{j-1})|}{f(a_{j-1})}.$$

- Approximation of $f''(\xi_{1,j})$ and $f''(\xi_{2,j})$: Define $\epsilon_{1,j}$ and $\epsilon_{2,j}$ by

$$\epsilon_{1,j} = f''(\xi_{1,j}) - f''(a_{j-1}), \text{ and } \epsilon_{2,j} = f''(\xi_{2,j}) - f''(a_j).$$

By uniform continuity of $f^{(2)} = f''$ on the compact set $[0, \tau]$, for every $\epsilon > 0$ there exists an $\eta = \eta_\epsilon > 0$ such that $|x - y| < \eta$ implies $|f''(x) - f''(y)| < \epsilon$. Fix $\epsilon > 0$ (to be chosen later). We have $\xi_{1,j}, \xi_{2,j} \in [a_{j-1}, a_{j+1}]$, where, by the proof of Lemma 4.1, if $k_n > 5\tilde{\gamma}_1 R$,

$$a_{j+1} - a_{j-1} = a_{j+1} - a_j + a_j - a_{j-1} \leq \frac{1}{k_n f(a_j^*)}(\sqrt{2} + 1)$$

$$\leq (\sqrt{2} + 1) \frac{1}{k_n f(\tau)}$$

$$\leq \frac{(\sqrt{2} + 1)R}{k_n}.$$

Thus, if we choose $k_n$ such that $k_n > \max\left(5\tilde{\gamma}_1 R, (\sqrt{2}+1)/\eta R\right)$, then $a_{j+1} - a_{j-1} < \eta$ for all $j = 1, \ldots, k$ and furthermore

$$\max\left\{|f''(\xi_{1,j}) - f''(a_{j-1})|, |f''(\xi_{2,j}) - f''(a_{j-1})|\right\} < \epsilon, \text{ for } j = 1, \ldots, k,$$

or, equivalently, $\max\{|\epsilon_{1,j}|, |\epsilon_{2,j}|\} < \epsilon, j = 1, \ldots, k$.

- Expanding $\Delta_{j+1} a$ around $\Delta_j a$: We have

$$\Delta_{j+1} a = a_{j+1} - a_j = a_j - a_{j-1} + [a_{j+1} - a_j - (a_j - a_{j-1})]$$

$$= \Delta_j a + \Delta_j a \left( \frac{a_{j+1} - a_j}{a_j - a_{j-1}} - 1 \right) = \Delta_j a + \Delta_j a \, \epsilon_{3,j}$$

where

$$\epsilon_{3,j} = \frac{a_{j+1} - a_j}{a_j - a_{j-1}} - 1 = \frac{f(a_j^*)}{f(a_{j+1}^*)} - 1 = \frac{f(a_j^*) - f(a_{j+1}^*)}{f(a_{j+1}^*)}$$
$$= \frac{-f'(a_j^{**})}{f(a_{j+1}^*)}(a_{j+1}^* - a_j^*).$$

Thus,

$$|\epsilon_{3,j}| \leq \frac{|f'(a_{j-1})|}{f(a_{j+1})} (a_{j+1} - a_{j-1})$$
$$= \frac{|f'(a_{j-1})|}{f(a_{j+1})} \left( \frac{1}{k_n f(a_j^*)} + \frac{1}{k_n f(a_{j+1}^*)} \right)$$
$$\leq 2 \frac{|f'(a_{j-1})|}{f^2(a_{j+1})} \frac{1}{k_n} = 2 \frac{|f'(a_{j-1})|}{f^2(a_{j-1})} \left( \frac{f(a_{j-1})}{f(a_{j+1})} \right)^2 \frac{1}{k_n}$$
$$\leq 2 \cdot 2^4 \frac{|f'(a_{j-1})|}{f^2(a_{j-1})} \frac{1}{k_n} = 32 \frac{|f'(a_{j-1})|}{f^2(a_{j-1})} \frac{1}{k_n} \leq 32\tilde{\gamma}_1 \frac{1}{k_n}.$$

Above, we have used the fact that $k_n > 5\tilde{\gamma}_1 R$ to be able to use the inequality $f(a_{j-1})/f(a_{j+1}) < 2^2$.
Now, expansion of $\beta_j$ yields, after straightforward algebra,

$$24\beta_j = \left[ -2M_j f''(a_{j-1})(\Delta_j a)^4 \right]$$
$$+ \left[ \epsilon_{1,j} \left( \frac{1}{2} + M_j \Delta_j a \right) (-\Delta_j a)^3 + \epsilon_{2,j} \left( \frac{1}{2} - M_j \Delta_j a \right) (\Delta_j a)^3 \right]$$
$$+ \left[ \left( \frac{1}{2} - M_j \Delta_j a \right) (3 + 3\epsilon_{3,j} + \epsilon_{3,j}^2)(f''(a_{j-1}) + \epsilon_{2,j}) \epsilon_{3,j} (\Delta_j a)^3 \right]$$
$$= T_{1,j} + T_{2,j} + T_{3,j}$$

where

$$\frac{|T_{1,j}|}{(\Delta_j a)^3} = 2|M_j| f''(a_{j-1})(\Delta_j a) \leq 2(\sqrt{2} + 1) \frac{|f'(a_{j-1})|}{f(a_{j-1})} f''(a_{j-1}) \frac{1}{k_n f(a_j^*)}$$
$$\leq 4(\sqrt{2} + 1) \frac{|f'(a_{j-1})|}{f^2(a_{j-1})} f''(a_{j-1}) \frac{1}{k_n}$$
$$\leq 4(\sqrt{2} + 1)\tilde{\gamma}_1 \bar{f}''_j \frac{1}{k_n} \leq 4(\sqrt{2} + 1)\tilde{\gamma}_1 \gamma_2 f(\tau)^3 \frac{1}{k_n}$$
$$\leq 2^{-1}(\sqrt{2} + 1)\tilde{\gamma}_1 \gamma_2 \tau^{-3} \equiv M_1 \frac{1}{k_n},$$

since $f(\tau) \leq (2\tau)^{-1}$ by (3.1), page 1669, Groeneboom, Jongbloed and Wellner [11],

$$\frac{|T_{2,j}|}{(\Delta_j a)^3} \leq 2 \left( \frac{1}{2} + \frac{2(\sqrt{2} + 1)\tilde{\gamma}_1}{k_n} \right) \epsilon \leq 2 \left( \frac{1}{2} + \frac{2(\sqrt{2} + 1)}{5R} \right) \epsilon = M_2 \epsilon,$$

and

$$\frac{|T_{3,j}|}{(\Delta_j a)^3} \leq \left(\frac{1}{2} + \frac{2(\sqrt{2}+1)\tilde{\gamma}_1}{k_n}\right)\left(3 + \frac{96\tilde{\gamma}_1}{k_n} + \frac{32^2\tilde{\gamma}_1^2}{k_n^2}\right)(\overline{f}''_j + \epsilon)\frac{1}{k_n}$$

$$\leq \left(\frac{1}{2} + \frac{2(\sqrt{2}+1)}{5R}\right)\left(3 + \frac{96}{5R} + \frac{32^2}{25R^2}\right) 2\gamma_2 f(\tau)^3 \frac{1}{k_n}$$

$$\leq \left(\frac{1}{2} + \frac{2(\sqrt{2}+1)}{5R}\right)\left(3 + \frac{96}{5R} + \frac{32^2}{25R^2}\right) 2^{-2}\gamma_2 \tau^{-3} \frac{1}{k_n} = M_3 \frac{1}{k_n}$$

if we choose $\epsilon < \gamma_2 f(\tau)^3 = \sup_{0 < t < \tau} f''(t)$ and again use $f(\tau) \leq (2\tau)^{-1}$. Note that by (21)

$$\frac{|t_j - r_j|}{(\Delta_j a)^3} \leq \frac{\max_{1 \leq i \leq k} |\mathcal{E}^{(1)}(a_i)|}{(\Delta_j a)^2}.$$

Thus, using (22) and combining the results obtained above, we can write for $j = 1, \ldots, k$,

$$\frac{|t_j - r_j|}{(\Delta_j a)^3} \leq \max_{1 \leq i \leq k-1} \frac{|\beta_i|}{(\Delta_j a)^2} \leq 24^{-1} \max_{1 \leq i \leq k-1} \frac{|T_{1,i}| + |T_{2,i}| + |T_{3,i}|}{(\Delta_i a)^3} \cdot \frac{|a|^3}{(\Delta_j a)^2}$$

$$\leq \left[(M_1 + M_3)\frac{1}{k_n} + M_2\, \epsilon\right] \frac{|a|^3}{(\Delta_j a)^2}$$

(23) $$= \left[(M_1 + M_3)\frac{1}{k_n} + M_2\, \epsilon\right] \frac{|a|^3}{(\Delta_j a)^3} \Delta_j a$$

But note that

$$\frac{|a|^3}{(\Delta_j a)^3} = \max_{1 \leq i \leq k}\left(\frac{\Delta_i a}{\Delta_j a}\right)^3 \leq \max_{1 \leq i \leq k}\left(\frac{f(a_j^*)}{f(a_i^*)}\right)^3 \leq \left(\frac{f(a_{j-1})}{f(\tau)}\right)^3$$

where

$$\frac{f(a_{j-1})}{f(\tau)} = \frac{f(a_{j-1})}{f(a_k)} = \frac{f(a_{j-1})}{f(a_j)} \cdot \frac{f(a_j)}{f(a_{j+1})} \cdots \frac{f(a_{k-1})}{f(a_k)}$$

and, for $l = 0, \ldots, k-1$,

$$\frac{f(a_l)}{f(a_{l+1})} = 1 + \frac{f(a_l) - f(a_{l+1})}{f(a_{l+1})}$$

$$= 1 + \frac{-f'(a_l^*)}{f(a_{l+1})}(a_{l+1} - a_l), \quad a_l^* \in [a_l, a_{l+1}]$$

$$= 1 + \frac{-f'(a_l^*)}{f(a_{l+1})f(a_l^{**})}\frac{1}{k_n}, \quad a_l^{**} \in [a_l, a_{l+1}]$$

$$\leq 1 + \frac{-f'(a_l)}{f(a_{l+1})f(a_l^{**})}\frac{1}{k_n}$$

$$= 1 + \frac{-f'(a_l)}{f^2(a_l)} \frac{f^2(a_l)}{f(a_{l+1})f(a_l^{**})} \frac{1}{k_n}$$

$$\leq 1 + \frac{\tilde{\gamma}_1}{4}\frac{1}{k_n}, \quad \text{if } k_n > 5\tilde{\gamma}_1 R.$$

Hence,

$$\frac{|a|^3}{(\Delta_j a)^3} \leq \left(1 + \frac{\tilde{\gamma}_1}{4}\frac{1}{k_n}\right)^{3(k_n+2-j)} \leq \left(1 + \frac{\tilde{\gamma}_1}{4}\frac{1}{k_n}\right)^{3(k_n+2)}$$

$$\leq \left(1 + \frac{\tilde{\gamma}_1}{4}\frac{1}{k_n}\right)^{3(k_n+2)} = \left(1 + \frac{\tilde{\gamma}_1}{4}\frac{1}{k_n}\right)^6 \left(1 + \frac{3\tilde{\gamma}_1}{4}\frac{1}{3k_n}\right)^{3k_n}$$

(24) $$\leq 2\left(1 + \frac{3\tilde{\gamma}_1}{4}\frac{1}{3k_n}\right)^{3k_n} \leq 2e^{3\tilde{\gamma}_1/4}$$

if $k_n \geq \tilde{\gamma}_1/(4(2^{1/6}-1))$ where we used $\log(1+x) \leq x$ for $x > 0$ in the last inequality. Combining (23) with (24), it follows that if we choose

$$k_n > \max\left\{5\tilde{\gamma}_1 R, \tilde{\gamma}_1/(4(2^{1/6}-1)), (\sqrt{2}+1)/\eta R\right\}$$

then

$$\frac{|t_j - r_j|}{(\Delta_j a)^3} \leq 4e^{3\tilde{\gamma}_1/4}\left[(M_1 + M_3)\frac{1}{k_n} + M_2\epsilon\right]\Delta_j a = o(\Delta_j a)$$

or

$$\frac{|t_j - r_j|}{(\Delta_j a)^4} = o(1)$$

where $o(1)$ is uniform in $j$. □

**Lemma 4.3.** *Under the hypotheses of Theorem 2.2,*

$$Pr\left(|T_j - t_j - (R_j - r_j)| \geq \delta_n p_n^3\right) \leq 4\exp\left(-\frac{(100)^{-1}n\delta_n^2 f^2(a_j^*)p_n^3}{1 + (1/30)p_n\delta_n f(a_j^*)}\right).$$

*Proof.* Write

$$W_j \equiv T_j - t_j - (R_j - r_j)$$
$$= -\left\{\frac{(\mathbb{Y}_n - Y)^{(1)}(a_{j-1}) + (\mathbb{Y}_n - Y)^{(1)}(a_j)}{2}\right.$$
$$\left. - \frac{(\mathcal{C}[\mathbb{Y}_n - Y])^{(1)}(a_{j-1}) + (\mathcal{C}[\mathbb{Y}_n - Y])^{(1)}(a_j)}{2}\right\}\Delta_j a$$
$$\equiv -\frac{1}{2}\left(\mathcal{E}^{(1)}_{\mathbb{Y}_n - Y}(a_{j-1}) + \mathcal{E}^{(1)}_{\mathbb{Y}_n - Y}(a_j)\right)\Delta_j a$$

where

$$\mathcal{E}^{(1)}_g(t) \equiv (g - \mathcal{C}[g])^{(1)}(t).$$

But for $g \in C^1[a_{j-1}, a_j]$ with $g^{(1)}$ of bounded variation,

$$g(t) = g(a_{j-1}) + g'(a_{j-1})(t - a_{j-1}) + \int_{a_{j-1}}^{t}(t-u)dg^{(1)}(u)$$
$$= P_j(t) + \int_{a_{j-1}}^{a_j} g_u(t)dg^{(1)}(u)$$

where $g_u(t) \equiv (t-u)_+ = (t-u)1_{[t \geq u]}$. Since $\mathcal{C}$ is linear and preserves linear functions

$$\mathcal{C}[g](t) = P_j(t) + \int_{a_{j-1}}^{a_j} \mathcal{C}g_u(t)dg^{(1)}(u),$$

and this yields

$$\mathcal{E}_g(t) = \int_{a_{j-1}}^{a_j} \mathcal{E}_{g_u}(t)dg^{(1)}(u)$$

and

$$\mathcal{E}_g^{(1)}(t) = \int_{a_{j-1}}^{a_j} \mathcal{E}_{g_u}^{(1)}(t)dg^{(1)}(u).$$

Applying this second formula to $g = \mathbb{Y}_n - Y$ yields the relation

$$\mathcal{E}_{\mathbb{Y}_n - Y}^{(1)}(t) = \int_{a_{j-1}}^{a_j} \mathcal{E}_{g_u}^{(1)}(t)d(\mathbb{F}_n - F)(u).$$

Now $g_u$ is absolutely continuous with $g_u(t) = \int_0^t g_u^{(1)}(s)ds$ where $g_u^{(1)}(t) = 1_{[t \geq u]}$, so by de Boor [5], (17) on page 56 (recalling that our $\mathcal{C} = I_4$ of de Boor),

$$\|\mathcal{E}_{g_u}^{(1)}\| = \|g_u^{(1)} - (\mathcal{C}[g_u])^{(1)}\|$$
$$\leq (19/4)\text{dist}(g_u^{(1)}, \$_3) \leq (19/4)\text{dist}(g_u^{(1)}, \$_2)$$
$$\leq (19/4)\omega(g_u^{(1)}, |a|) \leq (19/4) \leq 5.$$

Thus the functions $(u,t) \mapsto \mathcal{E}_{g_u}^{(1)}(t)\Delta_j a$ are bounded by a constant multiple of $\Delta_j a$, while the functions $h_{j,l}(u) = \mathcal{E}_{g_u}^{(1)}(a_l)1_{[a_{j-1},a_j]}(u)\Delta_j a$, $l \in \{j-1, j\}$ satisfy

$$Var[h_{j,l}(X)] \leq (\Delta_j a)^2 \int_{a_{j-1}}^{a_j} (19/4)^2 f(u)du \leq 5^2(\Delta_j a)^3 f(a_{j-1})$$
$$\leq 50p_n^3/f^2(a_j^*)$$

for $k \geq 5\tilde{\gamma}_1(F,\tau)R$ as in the proof of Lemma 3.1 in section 3. By applying Bernstein's inequality much as in the proof of Lemma 3.1 we find that

$$Pr\left(|\mathcal{E}_{\mathbb{Y}_n - Y}^{(1)}(a_l)| > \delta_n p_n^3\right)$$
$$\leq 2\exp\left(-\frac{n\delta_n^2 p_n^6/2}{50p_n^3/f(a_j^*)^2 + p_n(5/3)\delta_n p_n^3/f(a_j^*)}\right)$$
$$= 2\exp\left(-\frac{n\delta_n^2 f^2(a_j^*)p_n^3}{100 + (10/3)p_n f(a_j^*)\delta_n}\right)$$
$$= 2\exp\left(-\frac{(100)^{-1}n\delta_n^2 f^2(a_j^*)p_n^3}{1 + (1/30)p_n \delta_n f(a_j^*)}\right).$$

Thus it follows that

$$Pr\left(|W_j| > \delta_n p_n^3\right)$$
$$\leq Pr\left(|\mathcal{E}_{\mathbb{Y}_n - Y}^{(1)}(a_{j-1})| > \delta_n p_n^3\right)$$
$$+ Pr\left(|\mathcal{E}_{\mathbb{Y}_n - Y}^{(1)}(a_j)| > \delta_n p_n^3\right)$$
$$\leq 4\exp\left(-\frac{(100)^{-1}n\delta_n^2 f^2(a_j^*)p_n^3}{1 + (1/30)p_n \delta_n f(a_j^*)}\right).$$

This completes the proof of the claimed bound. □

**Lemma 4.4.** *Let $R(s,t)$ be defined by*

$$R(s,t) \equiv Ph_{s,t}$$
$$= \frac{1}{2}(F(t) + F(s))(t-s) - \int_s^t F(u)du, \qquad 0 \le s \le t < \infty.$$

*Then*

(25) $$R(s,t) \begin{cases} \le \frac{1}{12}f'(s)(t-s)^3 + \frac{1}{24}\sup_{s \le x \le t} f''(x)(t-s)^4 \\ \ge \frac{1}{12}f'(s)(t-s)^3 + \frac{1}{24}\inf_{s \le x \le t} f''(x)(t-s)^4. \end{cases}$$

**Remark.** It follows from the Hadamard-Hermite inequality that for $F$ concave, $R(s,t) \le 0$ for all $s \le t$; see e.g. Niculescu and Persson [19], pages 50 and 62-63 for an exposition and many interesting extensions and generalizations. Lemmas A4 and A5 give additional information under the added hypotheses that $F^{(2)}$ exists and $F^{(1)}$ is convex.

*Proof.* Since $g_s(t) \equiv R(s,t)$ has first three derivatives given by

$$g_s^{(1)}(t) = \frac{d}{dt}R_s(t) = \frac{1}{2}f(t)(t-s) + \frac{1}{2}(F(t) + F(s) - F(t))$$
$$= \frac{1}{2}f(t)(t-s) - \frac{1}{2}(F(t) - F(s)) \stackrel{t=s}{=} 0,$$
$$g_s^{(2)}(t) = \frac{d^2}{dt^2}R_s(t) = \frac{1}{2}f'(t)(t-s) + \frac{1}{2}(f(t) - f(t)) \stackrel{t=s}{=} 0,$$
$$g_s^{(3)}(t) = \frac{d^3}{dt^3}R_s(t) = \frac{1}{2}f''(t)(t-s) + \frac{1}{2}f'(t),$$

we can write $R(s,t)$ as a Taylor expansion with integral form of the remainder: for $s < t$,

$$R(s,t) = g_s(t) = g_s(s) + g_s'(s)(t-s) + \frac{1}{2!}g_s''(s)(t-s)^2$$
$$+ \frac{1}{2!}\int_s^t g_s^{(3)}(x)(t-x)^2 dx$$
$$= 0 + \frac{1}{2!}\int_s^t \left(\frac{1}{2}f''(x)(x-s) + \frac{1}{2}f'(x)\right)(t-x)^2 dx$$
$$= \frac{1}{4}\int_s^t f'(x)(t-x)^2 dx + \frac{1}{4}\int_s^t f''(x)(x-s)(t-x)^2 dx$$
$$= \frac{1}{4}\int_s^t \{f'(s) + f''(x^*)(x-s)\}(t-x)^2 dx$$
$$+ \frac{1}{4}\int_s^t f''(x)(x-s)(t-x)^2 dx$$
$$= \frac{1}{12}f'(s)(t-s)^3 + \frac{1}{4}\int_s^t \{f''(x^*) + f''(x)\}(x-s)(t-x)^2 dx$$

where $|x^* - x| \le |x - s|$ for each $x \in [s,t]$. Since $\int_s^t (x-s)(t-x)^2 dx = (t-s)^4/12$ we find that the inequalities (25) hold. □

**Lemma 4.5.** *Let $r_{n,i} \equiv P(h_{a_{i-1},a_i}) = R(a_{i-1}, a_i)$, $i = j, j+1$, $\underline{f}''_j = \inf_{t \in [a_{j-1},a_j]} f''(t)$ and $\overline{f}''_j = \sup_{t \in [a_{j-1},a_j]} f''(t)$. Then there exists $a_j^* \in [a_{j-1}, a_j] = I_j$ such that*

$$\frac{r_{n,j}}{(\Delta_j a)^3} - \frac{r_{n,j+1}}{(\Delta_{j+1} a)^3} \leq -\frac{1}{12} f''(a_j^*) \Delta_j a + \frac{1}{24}(\overline{f}''_j \Delta_j a - \underline{f}''_{j+1} \Delta_{j+1} a).$$

*Proof.* In view of (25), we have

$$r_{n,j} \begin{cases} \leq \frac{1}{12} f'(a_{j-1})(\Delta_j a)^3 + \frac{1}{24} \sup_{x \in I_j} f''(x)(\Delta_j a)^4 \\ \geq \frac{1}{12} f'(a_{j-1})(\Delta_j a)^3 + \frac{1}{24} \inf_{x \in I_j} f''(x)(\Delta_j a)^4, \end{cases}$$

$$r_{n,j+1} \begin{cases} \leq \frac{1}{12} f'(a_j)(\Delta_{j+1} a)^3 + \frac{1}{24} \sup_{x \in I_{j+1}} f''(x)(\Delta_{j+1} a)^4 \\ \geq \frac{1}{12} f'(a_j)(\Delta_{j+1} a)^3 + \frac{1}{24} \inf_{x \in I_{j+1}} f''(x)(\Delta_{j+1} a)^4, \end{cases}$$

and hence

$$\frac{r_{n,j}}{(\Delta_j a)^3} - \frac{r_{n,j+1}}{(\Delta_{j+1} a)^3}$$
$$\leq \frac{1}{12} f'(a_{j-1}) + \frac{1}{24} \sup_{x \in I_j} f''(x) \Delta_j a - \frac{1}{12} f'(a_j) - \frac{1}{24} \inf_{x \in I_{j+1}} f''(x) \Delta_{j+1} a$$
$$= -\frac{1}{12} f''(a_j^*) \Delta_j a + \frac{1}{24}(\overline{f}''_j \Delta_j a - \underline{f}''_{j+1} \Delta_{j+1} a), \text{ where } a_j^* \in I_j. \quad \square$$

## 5. Appendix 2: A "modernized" proof of Kiefer and Wolfowitz [14]

Define the following interpolated versions of $F$ and $\mathbb{F}_n$. For $k \geq 1$, let $a_j \equiv a_j^{(k)} \equiv F^{-1}(j/k)$ for $j = 1, \ldots, k-1$, and set $a_0 \equiv \alpha_0(F)$ and $a_k \equiv \alpha_1(F)$. Using the notation of de Boor [5], Chapter III, let $L^{(k)} = I_2 F$ be the piecewise linear and continuous function on $\mathbb{R}$ satisfying

$$L^{(k)}(a_j^{(k)}) = F(a_j^{(k)}), \qquad j = 0, \ldots, a_k.$$

Similarly, define $\mathbb{L}_n \equiv \mathbb{L}_n^{(k)} = I_2 \mathbb{F}_n$; thus

$$\mathbb{L}_n^{(k)}(x) = \mathbb{F}_n(a_j) + k\{\mathbb{F}_n(a_{j+1}) - \mathbb{F}_n(a_j)\}[L^{(k)}(x) - F(a_j)]$$

for $a_j \leq x \leq a_{j+1}$, $j = 0, \ldots, a_k$. We will eventually let $k = k_n$ and then write $p_n = 1/k_n$ (so that $F(a_{j+1}) - F(a_j) = 1/k_n = p_n$).

The following basic lemma due to Marshall [17] plays a key role in the proof.

**Lemma 5.1** (Marshall [17]). *Let $\Psi$ be convex on $[0,1]$, and let $\Phi$ be a continuous real-valued function on $[0,1]$. Let*

$$\overline{\Phi}(x) = \sup\{h(x): h \text{ is convex and } h(z) \leq \Phi(z) \text{ for all } z \in [0,1]\}.$$

*Then*

$$\sup_{0 \leq x \leq 1} |\overline{\Phi}(x) - \Psi(x)| \leq \sup_{0 \leq x \leq 1} |\Phi(x) - \Psi(x)|.$$

*Proof.* Note that for all $y \in [0,1]$, either $\overline{\Phi}(y) = \Phi(y)$, or $y$ is an interior point of a closed interval $I$ over which $\overline{\Phi}$ is linear. For such an interval, either $\sup_{x \in I} |\overline{\Phi}(x) -$ $\Psi(x)|$ is attained at an endpoint of $I$ (where $\overline{\Phi} = \Phi$), or it is attained at an interior point, where $\Psi < \overline{\Phi}$. Since $\overline{\Phi} \leq \Phi$ on $[0, 1]$, it follows that

$$\sup_{x \in I} |\overline{\Phi}(x) - \Psi(x)| \leq \sup_{x \in I} |\Phi(x) - \Psi(x)|.$$

Here is a second proof (due to Robertson, Wright and Dykstra [21], page 329) that does not use continuity of $\Phi$. Let $\epsilon \equiv \|\Phi - \Psi\|_\infty$. Then $\Psi - \epsilon$ is convex, and $\Psi(x) - \epsilon \leq \Phi(x)$ for all $x$. Thus for all $x$

$$\Phi(x) \geq \overline{\Phi}(x) \geq \Psi(x) - \epsilon,$$

and hence

$$\epsilon \geq \Phi(x) - \Psi(x) \geq \overline{\Phi}(x) - \Psi(x) \geq -\epsilon$$

for all $x$. This implies the claimed bound. □

**Main steps:**

A. By Marshall's lemma, for any concave function $h$, $\|\widehat{F}_n - h\| \leq \|\mathbb{F}_n - h\|$.
B. $P_F(A_n) \equiv P_F\{\mathbb{L}_n^{(k_n)}$ is concave on $[0, \infty)\} \nearrow 1$ as $n \to \infty$ if $k_n \equiv (C_0 \beta_1(F) \times n/\log n)^{1/3}$ for some absolute constant $C_0$.
C. On the event $A_n$, it follows from Marshall's lemma (step **A**) that

$$\begin{aligned}
\|\widehat{F}_n - \mathbb{F}_n\| &= \|\widehat{F}_n - \mathbb{L}_n^{(k_n)} + \mathbb{L}_n^{(k_n)} - \mathbb{F}_n\| \\
&\leq \|\mathbb{F}_n - \mathbb{L}_n^{(k_n)}\| + \|\mathbb{L}_n^{(k_n)} - \mathbb{F}_n\| \\
&= 2\|\mathbb{F}_n - \mathbb{L}_n^{(k_n)}\| \\
&= 2\|\mathbb{F}_n - \mathbb{L}_n^{(k_n)} - (F - L^{(k_n)}) + F - L^{(k_n)}\| \\
&\leq 2\|\mathbb{F}_n - \mathbb{L}_n^{(k_n)} - (F - L^{(k_n)})\| + 2\|F - L^{(k_n)}\| \\
&\equiv 2(D_n + E_n).
\end{aligned}$$

D. $D_n$ is handled by a standard "oscillation theorem"; $E_n$ is handled by an analytic (deterministic) argument.

*Proof of (1) assuming B holds.* Using the notation of de Boor [5], chapter III, we have

$$\mathbb{F}_n - F - (\mathbb{L}_n - L) = \mathbb{F}_n - F - I_2(\mathbb{F}_n - F).$$

But by (18) of de Boor [5], page 36, $\|g - I_2 g\| \leq \omega(g; |a|)$ where $\omega(g; |a|)$ is the oscillation modulus of $g$ with maximum comparison distance $|a| = \max_j \Delta a_j$ (and note that de Boor's proof does not involve continuity of $g$). Thus it follows immediately that

$$\begin{aligned}
D_n &\equiv \|\mathbb{F}_n - F - (\mathbb{L}_n - L)\| \\
&= \|\mathbb{F}_n - F - I_2(\mathbb{F}_n - F)\| \\
&\leq \omega(\mathbb{F}_n - F; |a|) \stackrel{d}{=} n^{-1/2} \omega(\mathbb{U}_n; p_n)
\end{aligned}$$

where $\mathbb{U}_n \equiv \sqrt{n}(\mathbb{G}_n - I)$ is the empirical process of $n$ i.i.d. Uniform$(0, 1)$ random variables. From Stute's theorem (see e.g. Shorack and Wellner [22], Theorem 14.2.1, page 542), $\limsup \omega(\mathbb{U}_n; p_n)/\sqrt{2 p_n \log(1/p_n)} = 1$ almost surely if $p_n \to 0$, $n p_n \to \infty$ and $\log(1/p_n)/n p_n \to 0$. Thus we conclude that

$$\|\mathbb{F}_n - F - (\mathbb{L}_n - L)\| = O(n^{-1/2} \sqrt{p_n \log(1/p_n)}) = O((n^{-1} \log n)^{2/3})$$

almost surely as claimed.

To handle $E_n$, we use the bound given by de Boor [5], page 31, (2): $\|g - I_2 g\| \leq 8^{-1}|a|^2 \|g''\|$. Applying this to $g = F$, $I_2 g = L^{(k)}$ yields

$$\|F - L^{(k)}\| = \|F - I_2 F\| \leq \frac{1}{8}|a|^2 \|F''\|$$
$$\leq \frac{1}{8}\gamma_1(F)p_n^2 = O((n^{-1} \log n)^{2/3}).$$

Combining the results for $D_n$ and $E_n$ yields the stated conclusion. $\square$

It remains to show that **B** holds. To do this we use the following lemma.

**Lemma 5.2.** *If $p_n \to 0$ and $\delta_n \to 0$, then for the uniform$(0,1)$ d.f. $F = I$,*

$$P(|\mathbb{G}_n(p_n) - p_n| \geq \delta_n p_n) \leq 2 \exp(-\frac{1}{2}np_n\delta_n^2(1+o(1)))$$

*where the $o(1)$ term depends only on $\delta_n$.*

*Proof.* From Shorack and Wellner [22], Lemma 10.3.2, page 415,

$$P(\mathbb{G}_n(p_n)/p_n \geq \lambda) \leq P\left(\sup_{p_n \leq t \leq 1} \frac{\mathbb{G}_n(t)}{t} \geq \lambda\right) \leq \exp(-np_n h(\lambda))$$

where $h(x) = x(\log x - 1) + 1$. Hence

$$P\left(\frac{\mathbb{G}_n(p_n) - p_n}{p_n} \geq \lambda\right) \leq \exp(-np_n h(1+\lambda))$$

where $h(1+\lambda) \sim \lambda^2/2$ as $\lambda \downarrow 0$, by Shorack and Wellner [22], (11.1.7), page 44. Similarly, using Shorack and Wellner [22], (10.3.6) on page 416,

$$P\left(\frac{p_n - \mathbb{G}_n(p_n)}{p_n} \geq \lambda\right) = P\left(\frac{p_n}{\mathbb{G}_n(p_n)} \geq \frac{1}{1-\lambda}\right) \leq \exp(-np_n h(1-\lambda))$$

where $h(1-\lambda) \sim \lambda^2/2$ as $\lambda \searrow 0$. Thus the conclusion follows with $o(1)$ depending only on $\delta_n$. $\square$

Here is the lemma which is used to prove **B**.

**Lemma 5.3.** *If $\beta_1(F) > 0$ and $\gamma_1(F) < \infty$, then for $k_n$ large,*

$$1 - P(A_n) \leq 2k_n \exp(-n\beta_1^2(F)/80k_n^3).$$

*Proof.* For $1 \leq j \leq k_n$, write

$$T_{n,j} \equiv \mathbb{F}_n(a_j) - \mathbb{F}_n(a_{j-1}), \qquad \Delta_j a \equiv a_j - a_{j-1}.$$

By linearity of $L_n^{(k_n)}$ on the sub-intervals $[a_{j-1}, a_j]$,

$$A_n = \bigcap_{j=1}^{k_n-1} \left\{\frac{T_{n,j}}{\Delta_j a} \geq \frac{T_{n,j+1}}{\Delta_{j+1} a}\right\} \equiv \bigcap_{j=1}^{k_n-1} B_{n,j}.$$

Suppose that

(26) $\quad |T_{n,i} - 1/k_n| \leq \delta_n/k_n, \quad i = j, j+1;\quad$ and $\quad \dfrac{\Delta_{j+1}a}{\Delta_j a} \geq 1 + 3\delta_n.$

Then
$$T_{n,j} \geq \frac{1}{k_n} - \frac{\delta_n}{k_n} = \frac{1-\delta_n}{k_n}, \qquad T_{n,j+1} \leq \frac{1+\delta_n}{k_n},$$

and it follows that for $\delta_n \leq 1/3$

$$T_{n,j}\frac{\Delta_{j+1}a}{\Delta_j a} \geq \frac{1-\delta_n}{k_n}(1+3\delta_n) \geq \frac{1-\delta_n}{k_n}\frac{1+\delta_n}{1-\delta_n} \geq T_{n,j+1}.$$

$[1+3\delta \geq (1+\delta)/(1-\delta)$ iff $(1+2\delta-3\delta^2) \geq 1+\delta$ iff $\delta - 3\delta^2 \geq 0$ iff $1 - 3\delta \geq 0.]$
Now the $\Delta$ part of (26) holds for $1 \leq j \leq k_n - 1$ provided $\delta_n \leq \beta_1(F)/6k_n < 1/3$.
Proof: Since

$$\frac{d}{dt}F^{-1}(t) = \frac{1}{f(F^{-1}(t))} \quad \text{and} \quad \frac{d^2}{dt^2}F^{-1}(t) = -\frac{f'}{f^3}(F^{-1}(t))$$

we can write

$$\Delta_{j+1}a = F^{-1}(\frac{j+1}{k}) - F^{-1}(\frac{j}{k}) = k_n^{-1}\frac{1}{f(a_j)} + \frac{1}{2k_n^2}\left(\frac{-f'(\xi)}{f^3(\xi)}\right)$$

for some $a_j \leq \xi \leq a_{j+1}$, and

$$\Delta_j a \leq k_n^{-1}\frac{1}{f(a_j)}.$$

Combining these two inequalities yields

$$\begin{aligned}\frac{\Delta_{j+1}a}{\Delta_j a} &\geq 1 + (2k_n)^{-1}f(a_j)\left(\frac{-f'(\xi)}{f^3(\xi)}\right) \\ &\geq 1 + \frac{1}{2k_n}\left(\frac{-f'(\xi)}{f^2(\xi)}\right) \geq 1 + \frac{1}{2k_n}\beta_1(F) \\ &= 1 + 3\delta_n\end{aligned}$$

if $\delta_n \equiv \beta_1(F)/(6k_n)$.

Thus we conclude that

$$\begin{aligned}1 - P(A_n) = P(\bigcup_{j=1}^{k_n-1} B_{n,j}^c) &\leq \sum_{j=1}^{k_n-1} P(B_{n,j}^c) \\ &\leq \sum_{j=1}^{k_n-1} 2P(|T_{n,j} - 1/k_n| > \delta_n/k_n) \\ &\leq k_n 4\exp(-2^{-1}np_n\delta_n^2 1 + o(1))) = 4k_n\exp(-n\beta_1^2(F)/80k_n^3).\end{aligned}$$

by using Lemma 5.2 and for $k_n$ sufficiently large (so that $(1 + o(1)) \geq 72/80$). $\square$

Putting these results together yields Theorem 1.1.

**Acknowledgments.** The second author owes thanks to Lutz Dümbgen and Kaspar Rufibach for the collaboration leading to the analogue of Marshall's lemma which is crucial for the development here. The second author thanks Piet Groeneboom for many stimulating discussions about estimation under shape constraints, and particularly about estimation of convex densities, over the past fifteen years.